\documentclass[journal ]{new-aiaa}
\usepackage[utf8]{inputenc}
\usepackage{textcomp}

\usepackage{graphicx}
\usepackage{amsmath}
\usepackage[version=4]{mhchem}
\usepackage{siunitx}
\usepackage{longtable,tabularx}
\usepackage{soul}
\setstcolor{red}

\usepackage{algorithm}
\usepackage{algpseudocode}
\usepackage{bm}
\usepackage{float}
\usepackage{subcaption}
\usepackage{hhline}
\usepackage{upgreek}
\usepackage{makecell}
\usepackage{hyperref}
\hypersetup{colorlinks=True,urlcolor=blue,linkcolor=blue,citecolor=blue}

\usepackage{anyfontsize}
\usepackage{datetime}
\usepackage[firstpageonly]{draftwatermark}
\SetWatermarkColor{red}
\SetWatermarkText{Submitted for consideration to Journal of Guidance, Control, and Dynamics, draft: \today ~ \currenttime}
\SetWatermarkAngle{0}
\SetWatermarkFontSize{5pt}
\SetWatermarkVerCenter{0.97\paperheight}

\usepackage{lineno}

\title{Advancing Averaged Primer Vector Theory \\with Bang-Bang Control and Eclipsing}

\author{Noah Lifset\footnote{Corresponding Author; Graduate Student, Aerospace Engineering and Engineering Mechanics, \href{mailto:noahlifset@utexas.edu} {\texttt{noahlifset@utexas.edu}}, \url{https://orcid.org/0000-0003-3397-7021}}, Ryan P. Russell\footnote{Professor, Aerospace Engineering and Engineering Mechanics, \href{mailto:ryan.russell@austin.utexas.edu} {\texttt{ryan.russell@austin.utexas.edu}}, \url{http://orcid.org/0000-0001-7672-0408}}}

\affil{University of Texas at Austin, Austin, Texas, 78712}

\begin{document}

\maketitle
\begin{abstract}
Primer vector theory using averaged dynamics is well suited for optimizing low-thrust, many-revolution spacecraft trajectories, but is difficult to implement in a way that maintains both optimality and computational efficiency. An improved model is presented that combines advances from several past works into a general and practical formulation for minimum-fuel, perturbed Keplerian dynamics. The model maintains computational efficiency of dynamics averaging with optimal handling of the eclipsing constraint and bang-bang control through the use of the Leibniz integral rule for multi-arc averaging. A subtle, but important singularity arising from the averaged eclipsing constraint is identified and fixed. A maximum number of six switching function roots per revolution is established within the averaged dynamics. This new theoretical insight provides a practical upper-bound on the number of thrusting arcs required for any low-thrust optimization problem. Variational equations are provided for fast and accurate calculation of the state transition matrix for use in targeting and optimization. The dynamics include generic two-body perturbations and an expanded state to allow for sensitivity calculations with respect to launch date and flight time. The new model is illustrated on a GTO to GEO transfer, including up to 486 revolutions.
\end{abstract}

\section{Introduction}

The increasing use of low-thrust, high-efficiency propulsion systems is making long flight-time, many-revolution trajectories (on the order of hundreds of revolutions) more relevant.\footnotemark[3]  \footnotetext{\footnotemark[3]A preliminary version of this paper was presented at the Jan 2024 AIAA SciTech Forum in Orlando, FL (paper number AIAA 2024-1279 \cite{lifset})} While it can often be sufficient to use a simple control law for these many-revolution trajectories \cite{petropoulos2004low, PPE, wu2021rapid}, optimization can provide significant savings in fuel and/or time of flight, and potentially enable new missions or extend the capability of existing missions. In the mission design process, the need for repeated trajectory optimization (e.g. trade studies) requires  optimization methods that are efficient and robust.

The indirect approach using primer vector theory \cite{lawden_1963} is a popular choice for handling low-thrust, many revolution spacecraft trajectories \cite{surveypaper}. This method transforms the optimization problem into a nonlinear two-point boundary value problem by introducing Lagrange multipliers (or costates) to solve for the optimal control law throughout the trajectory. The problem is reduced to finding the initial values of the costates. While the indirect method benefits from low dimensionality, the problems are notoriously hard due to high numerical sensitivity and complexity associated to incorporating constraints. In 1965 and 1967, Edelbaum used optimal control to analytically solve for low-thrust transfers between nearby elliptical orbits \cite{edelbaum1965,edelbaum1967general}. His approach was extended to the general case of transfers between arbitrary elliptical orbits by Marec and Vinh \cite{marec1977optimal}. Kechichian improved on these methods by using equinoctial elements to reduce the number of singularities \cite{kechichian1997trajectory}. More recent works have successfully applied the indirect method to many revolution trajectories using a variety of approaches, such as optimal handling of bang-bang control with switching detection \cite{wangtopputo,russell2007primer}, a smooth approximation of bang-bang control \cite{TAHERI2020151,taheri}, and interior point constraints \cite{wangtopputo,PONTANI2022549,cerf2018}.

In 1973, Edelbaum et al. introduced the use of dynamics averaging in optimal control for minimum time problems \cite{edelbaum}, which increased optimization speed and robustness. This averaging method was then used in the SEPSPOT software developed by him, Sackett, and Malchow \cite{sackett1975solar}. Marec and Vinh extended optimal control with dynamics averaging to the minimum fuel problem \cite{marec}. Kluever \cite{kluever} builds on SESPOT with a model that calculates the control profile by discretizing the costate time histories and solving the resulting nonlinear problem, although the model does not include variations in right ascension of the ascending node or argument of periapse. Geffroy and Epenoy \cite{geffroy} incorporated constraints on thrust direction and placed optimal control with dynamics averaging on more sound mathematical footing using Chaplais' theory of averaging \cite{chaplais}. While the effect of Earth's shadow had been included in previous averaging models, it was incorporated in a technically optimal manner for the first time by Dargent \cite{dargent} using a smoothed model for the constraint (Dargent also investigates the transformation between averaged and unaveraged solutions in \cite{dargent2}). Mazzini \cite{mazzini,mazzini2} explored the effect of the eclipsing constraint in averaging in more detail using a discontinuous model for the constraint. Tarzi \cite{tarzi} incorporated solar radiation pressure into the model and introduced a multi-arc approach to averaging for bang-bang control, but did not account for the Earth's shadow in an optimal way. Wu et al. \cite{wu_avg} develop a method to analytically calculate the roots of the switching function in the averaging context and use a neural network for initial guess generation, but does not include eclipsing in an optimal manner or in examples. Many of these more recent averaging methods are compared by Petropolous \cite{petropoulos_averaging}. The objective of this current work is to combine and extend the advances for the averaged dynamics case into a general and practical formulation for minimum fuel problems. The new model includes optimal handling of the eclipsing constraint, multi-arc averaging with novel application of the Leibniz integral rule, inclusion of generic dynamics perturbations, and variational equations for averaged dynamics for the first time.

Dynamics averaging removes the fast periodic changes in the states and costates, leaving only the slowly changing secular motion. The resulting averaged trajectory will approximate the unaveraged trajectory to within $\epsilon$, for integration times of up to $1/\epsilon$ \cite{sanders2007averaging}, where $\epsilon$ is typically a small parameter defined in terms of perturbation amplitude (control or dynamics). A numerical integrator can then take much larger time steps, which both reduces the numerical sensitivity to the initial costates and reduces computation time to solve the resulting system of ordinary differential equations. Dynamics averaging has also been used effectively in collocation methods that discretize and optimize the states and costates directly. Falck \cite{falck2012} builds on the averaging work done by Kluever \cite{kluever} to develop a general collocation method for the averaged minimum time problem. Olikara \cite{olikara2018,olikara2019} was the first to introduce a collocation method for the averaged minimum fuel problem. Shannon \cite{shannon2024} improves on Olikara's method for optimal coasting in collocation and introduces a hybrid approach that switches to unaveraged dynamics in regimes where the averaging approximation deteriorates. The collocation approach allows for straightforward inclusion of path constraints, such as eclipsing, but typically requires the number of control variables to scale with time of flight.

In practice, dynamics averaging is done by taking the Hamiltonian (with the embedded primer vector feedback law) and integrating over one period in time to find the averaged Hamiltonian. Derivatives needed for equations of motion (and the variational equations) are then calculated using the averaged Hamiltonian. Following the methods of Kechichian \cite{kechichian} and Edelbaum et al. \cite{edelbaum}, the averaging integrals are approximated using a weighted sum with Gaussian Quadrature because they rarely can be calculated analytically. For the minimum fuel problem with bang-bang control, a significant number of quadrature points are required to accurately capture the thrusting and coasting proportions of the averaging period. The multi-arc averaging approach introduced by Tarzi \cite{tarzi} for bang-bang control separates the integral over the Hamiltonian into different thrusting or ballistic ``arcs''. The multi-arc approach reduces the number of quadrature points required to achieve a specific level of accuracy at the cost of root finding in each dynamics function call. 

One of the most important constraints for low-thrust trajectory design is the effect of the Earth's shadow, which typically requires a spacecraft with an electric propulsion engine to stop thrusting. Incorporating this constraint in a way that maintains optimality is not simple in general, and especially difficult when using dynamics averaging. In a proper optimal control derivation, the eclipsing constraint is directly incorporated into the dynamics so that the control law is optimal with respect to the constraint \cite{singh_eclipse,WU2021107146}. If, instead, the constraint is simply enforced in flight, using a thrust profile that did not consider eclipsing, then the control law is technically no longer optimal. Mazzini describes the importance of the dynamics-based control constraint with multiple examples, showing how it can provide significant fuel savings depending on the Sun position at launch \cite{mazzini}. When the constraint is properly incorporated, the costates exhibit a discrete ``jump'' at eclipse entrance and exit. In averaged dynamics, these discrete jumps in a costate variable are not easily incorporated. Mazzini \cite{mazzini} develops a method to directly calculate the discrete jumps by integrating approximate costate dynamics through a transitionary layer for eclipsing. The averaged costate jumps are then calculated by applying the averaging integral to costate equations of motion that include a Dirac delta function. The resulting form of the averaged costate jumps is sufficient for propagation of an averaged trajectory, but does not address calculation of the variational equations. The model introduced in this current work incorporates the costate jumps into the averaged dynamics using the Leibniz integral rule applied to multi-arc averaging. This new approach uses no approximations and enables straightforward calculation of the variational equations.

In the current work, the model includes the multi-arc averaging approach introduced by Tarzi \cite{tarzi} along with a novel application of the Leibniz integral rule to automatically incorporate the averaged costate jumps. The new approach is simple, efficient, optimal with respect to the eclipsing constraint, and allows for straight-forward calculation of the variational equations. The averaged costate jumps, however, approach a singularity as eclipsing arcs approach zero length. This singularity, also identified in Mazzini's model, can significantly slow down or prevent propagation with a variable step integrator and result in accumulation of error. The new averaging model removes this singularity through a simple redefinition of the eclipsing constraint.
 
The focus of the current work is on the model for dynamics and sensitivities of averaged low-thrust trajectories using dynamics with eclipsing and perturbations. For the first time in the literature, variational equations of the full augmented state of the averaged dynamics are explicitly provided. An example of an optimal transfer from geostationary transfer orbit (GTO) to geostationary orbit (GEO) is shown as a representative application of the model and to validate the accuracy of the averaged equations. This validation is accomplished through comparison with an unaveraged trajectory propagated with approximate bang-bang control and eclipsing. An example of a 486-revolution optimal transfer from GTO to GEO with lower thrust and longer time-of-flight is then shown to showcase the capabilities of the averaged dynamics. The remainder of the paper is organized as follows. First, in section \ref{sec:optimal_control} the optimal control derivation is shown for unaveraged dynamics with bang-bang control and the eclipsing constraint. In section \ref{sec:averaging}, the averaged dynamics are introduced, multi-arc averaging is explained, and a new switching function polynomial is derived. In section \ref{sec:singularity}, the eclipsing singularity is explained and fixed. In section \ref{sec:variational}, the variational equations are calculated for the averaged dynamics. Section \ref{sec:results} presents the GTO to GEO examples. Approximate unaveraged dynamics are introduced in Section \ref{sec:approx_unaveraged} solely for ease of comparison to the averaged model. Concluding remarks are then given in section \ref{sec:conclusion}.

\section{Optimal Control} \label{sec:optimal_control}

Similar to prior works on this many-revolution low-thrust problem \cite{junkinsTaheri,taheri,graham2016minimum}, modified equinoctial orbital elements are used to reduce the number of singularities compared to Keplerian orbital elements \cite{walker1986set}. These elements are defined in terms of the Keplerian orbital elements as:
\begin{align}
\begin{split}
    p &= a[1-e^2]\\
    f &= e \cos(\omega + \Omega)\\
    g &= e \sin(\omega + \Omega)\\
    h &= \tan(i/2) \cos(\Omega)\\
    k &= \tan(i/2) \sin(\Omega)\\
    L &= \nu + \omega + \Omega
\end{split}
\end{align}
where \textit{a} is the semimajor axis, \textit{i} is the inclination, \textit{e} is the eccentricity, $\omega$ is the argument of perigee (AOP), $\Omega$ is the right ascension of the ascxcending node (RAAN), and $\nu$ is the true anomaly (TA). Current time, $t$, and overall time of flight, $\alpha$, are included in the state variables:
\begin{align}
    \text{\textbf{x}} &= [p, f, g, h, k, L, t, \alpha]^{\intercal}
\end{align}

The expanded state vector allows for better handling of boundary conditions and optimization, because all variables of interest in a given optimization problem are part of the state. While $t$ is linear in time and $\alpha$ is constant, their explicit inclusion simplifies addition of time dependent perturbations and calculation of the variational equations that include time-based sensitivities. A new time-like variable, $\tau$, is introduced that goes from 0 to 1 over the course of the trajectory with derivatives calculated per the chain rule:
\begin{align}
    \tau &= t \ / \ \alpha \\
    \frac{d}{d \tau} &= \frac{d t}{d \tau}\frac{d}{d t} = \alpha \frac{d}{d t}
\end{align}

\begin{figure}[h]
    \centering
    \vspace{60pt}
    \includegraphics[scale=0.45,trim={5cm 5cm 5cm 5cm}]{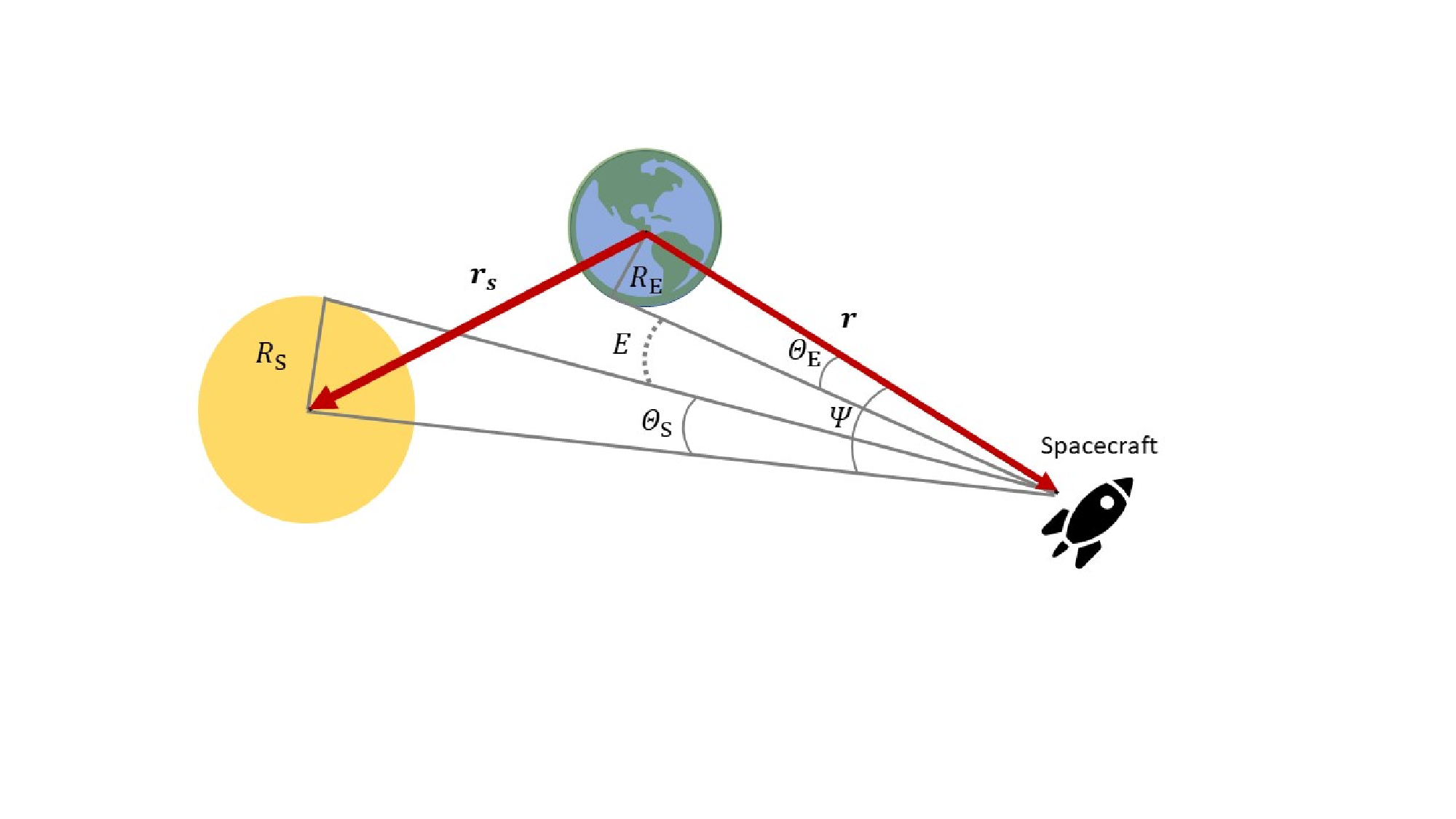}
    \caption{Angles used to calculate eclipsing function in conical shadow model.}
    \label{fig:shadow}
\end{figure}

In order to incorporate the effects of the Earth's shadow directly into the dynamics, the conventional geometric ``eclipsing function'' is used \cite{singh_eclipse}. The eclipsing function is a smooth function of the geocentric satellite position, the geocentric Sun position, the radius of the Sun, the radius of the Earth, and indirectly the state variable time. A conical shadow model is used, where the engine fully shuts off in umbra, penumbra, and antumbra (see Fig. \ref{fig:shadow})\cite{williamsEclipse}. The eclipsing function is written as
\begin{align}
\begin{split}
    &E(\text{\textbf{r}}_{\text{S}}(t),\text{\textbf{r}}(\textbf{x})) = \mathit{\Theta}_{\text{S}} + \mathit{\Theta}_{\text{E}} - \Psi \\
    &\mathit{\Theta}_{\text{S}} = \arcsin\left(\frac{R_{\text{S}}}{||\bm{r}_{\text{S}}(t)-\bm{r}(\textbf{x})||}\right) \\
    &\mathit{\Theta}_{\text{E}} = \arcsin\left(\frac{R_{\text{E}}}{r(x)}\right) \\
    &\mathit{\Psi} = \arccos\left(-\bm{\hat{r}}(\textbf{x}) \cdot \left[\frac{\bm{r}_{\text{S}}(t)-\bm{r}(\textbf{x})}{||\bm{r}_{\text{S}}(t)-\bm{r}(\textbf{x})||}\right]\right) \label{eq:eclipsing_func1}
\end{split}
\end{align}

\noindent where $\bm{r}_{\text{S}}$ is the geocentric Sun position, $\bm{r}$ is the geocentric satellite position, $R_{\text{S}}$ is the Sun radius, and $R_{\text{E}}$ is the Earth radius. Note that all equations in the paper use parentheses $``()"$ for functions, and square brackets $``[]"$ and curly brackets $``\{\}"$ for algebraic groupings.

A new function, $k_e$(\textit{E}), is defined that equals either 0 or 1 depending on the sign of the eclipsing function:
\begin{align}
    k_e = 
    \begin{cases}
        1 & \text{if } E < 0 \\
        0 & \text{if } E > 0
    \end{cases} \label{eq:ke}
\end{align}

The spacecraft thrust is defined in the radial-transverse-normal (RTN) reference frame. This frame is defined by the spacecraft radial vector, the transverse velocity vector, and the normal vector that completes the orthogonal basis. The controllable variables are the thrust direction, $\Hat{\bm{u}}$, and a thrust modulating parameter, $\sigma \in [0,1]$, with maximum and minimum thrusts given as $T_{\text{max}}$ and $T_{\text{min}}$, respectively. The $k_e$ function is incorporated directly into the dynamics. The equations of motion are defined:
\begin{align}
    &\frac{d m}{d \tau} = -\alpha \frac{T}{c}, \ \ c = g_0 \ I_{\text{SP}} \label{eq:unavg_mass_dynamics}\\
    &\frac{d \text{\textbf{x}}}{d \tau} = \alpha \text{\textbf{A}} + \alpha \text{\textbf{B}} \left[\hat{\bm{u}} \frac{T}{m} + \bm{\gamma}\right] \label{eq:unavg_state_dynamics}\\
    &T = T_{\text{min}} + [T_{\text{max}} - T_{\text{min}}] k_e \sigma 
    \label{eq:unavg_thrust_dynamics}
\end{align}

\noindent where $\bm{\gamma}$ is an arbitrary perturbation, $g_0$ is the gravitation acceleration at Earth's surface, and $I_{\text{SP}}$ is the specific impulse. The matrices \textbf{A}, \textbf{B} contain the state dynamics equations for modified equinoctial orbital elements written in the RTN frame. They are derived by Kechichian \cite{kechichian} and shown in the appendix. The arbitrary perturbation, $\bm{\gamma}$, and the thrust unit vector, $\hat{\bm{u}}$, are defined in the RTN frame, and $\bm{\gamma}$ must be a smooth function of state variables. This perturbation must be small enough to not significantly change the mean orbital elements over one period for the averaging assumptions to hold. Examples of generic perturbations that can be included are oblateness terms for the primary body, simple solar radiation pressure models, or simple drag models. In later examples, the $J_2$ perturbation is included.

The minimum fuel cost function \cite{bertrandepenoy} used in this model is defined as:
\begin{align}
\begin{split}
    J &= \int_0^1 \alpha \frac{T}{c} d\tau = \int_0^1 \left\{\alpha \frac{T_{\text{min}}}{c} + \alpha \frac{T_{\text{max}}-T_{\text{min}}}{c}k_e\sigma \right\}d\tau
    \label{eq:PI}
\end{split}
\end{align}

Following Eqs. (\ref{eq:unavg_state_dynamics}) and (\ref{eq:PI}), the Hamiltonian is therefore defined as:
\begin{align}
    H &= \alpha \frac{T}{c} + \alpha \bm{\uplambda}^{T} \text{\textbf{A}} + \alpha \bm{\uplambda}^{T} \text{\textbf{B}} \left[\hat{\bm{u}}  \frac{T}{m} +\gamma\right] - \alpha \lambda_m \frac{T}{c} \label{eq:hamiltonian_unavg}
\end{align}

\noindent where $\bm{\uplambda}$ is the costate vector for the 6 orbital elements and two time related variables, and $\lambda_m$ is the costate for the mass. The costate dynamics are then calculated as:
\begin{align}
    \frac{d \bm{\uplambda}}{d \tau} &= -\frac{\partial H}{\partial \text{\textbf{x}}}, \ \ \ \frac{d \bm{\uplambda}_m}{d \tau} = -\frac{\partial H}{\partial m} \label{eq:unavg_costate_dynamics}
\end{align}

To derive the optimal control law, the Hamiltonian is minimized with respect to the controls, thrust direction and thrust magnitude \cite{kirk1970optimal}. For the direction, the thrust vector is anti-aligned with $\bm{\uplambda}^{\intercal}$\textbf{B} \cite{lawden_1963}:
\begin{align}
    \hat{\bm{u}} &= - \frac{\bm{\uplambda}^{\intercal} \text{\textbf{B}}}{||\bm{\uplambda}^{\intercal} \text{\textbf{B}} ||}
\end{align}

\noindent The so-called "primer vector" \cite{lawden_1963} is substituted back into Eq. (\ref{eq:middle_H}) for an updated Hamiltonian:
\begin{align}
\begin{split}
    H =& \alpha \frac{T}{c} + \alpha \bm{\uplambda}^{T} \text{\textbf{A}} - \alpha ||\bm{\uplambda}^{T} \text{\textbf{B}}|| \frac{T}{m} +\alpha \bm{\uplambda}^{T} \text{\textbf{B}} \gamma - \alpha \lambda_m \frac{T}{c} \\
    =&\alpha \frac{T_{\text{min}}}{c} + \alpha \bm{\uplambda}^{T} \text{\textbf{A}} - \alpha ||\bm{\uplambda}^{T} \text{\textbf{B}}|| \frac{T_{\text{min}}}{m} +\alpha \bm{\uplambda}^{T} \text{\textbf{B}} \gamma - \alpha \lambda_m \frac{T_{\text{min}}}{c} \\
    &+ \alpha k_e [T_{\text{max}} - T_{\text{min}}] \frac{1}{c} \left[-\frac{c}{m} ||\bm{\uplambda}^{T} \text{\textbf{B}}|| - \lambda_m + 1 \right] \sigma \label{eq:middle_H}
\end{split}
\end{align}

The optimal $\sigma$ is determined by the sign of the quantity multiplying it, the so-called "switching function":
\begin{align}
    S' &= \alpha k_e [T_{\text{max}} - T_{\text{min}}] \frac{1}{c} \left[-\frac{c}{m} ||\bm{\uplambda}^{T} \text{\textbf{B}}|| - \lambda_m + 1 \right]
    \label{eq:switch_full}
\end{align}

\noindent This equation for the switching function can be further simplified. The quantities $\alpha, [T_{\text{max}} - T_{\text{min}}],$ and $c$ are all greater than zero and will not affect the sign of Eq. (\ref{eq:switch_full}). The variable $k_e$ will either be 1 and not affect the sign of Eq. (\ref{eq:switch_full}) or it will be 0, in which case $\sigma$ will have no effect on the equations of motion for the states or the costates. The latter case is technically a singular arc, but does not require further investigation because the singular control has no effect on the equations of motion. The optimal $\sigma$ is thus determined by a simplified switching function:
\begin{align}
    S &= -\frac{c}{m} ||\bm{\uplambda}^{T} \text{\textbf{B}}|| - \lambda_m + 1
    \label{eq:switch}
\end{align}
which leads to so-called ``bang-bang'' control:
\begin{align}
    \sigma &= \begin{cases}
        1 & \text{if } S < 0 \text{ and } k_e=1\\
        0 & \text{if } S > 0 \text{ and } k_e=1\\
        N/A & \text{if } k_e=0
    \end{cases} \label{eq:sigma}
\end{align}

\noindent This optimal $\sigma$ is then substituted back into the Hamiltonian in Eq. (\ref{eq:middle_H}). The derived optimal control law is now optimal with respect to the eclipsing constraint.

The traditional approach to maintaining optimality with a control constraint like eclipsing is to adjoin the Hamiltonian with a constraint function \cite{brysonho,hull,clarke2013}. The adjoined Hamiltonian is then used to calculate the Weierstrass-Erdmann corner conditions, which determine how the costates ``jump'' at eclipse entry and exit:
\begin{align}
    \bm{\uplambda}^T(\tau^-) = \bm{\uplambda}^T(\tau^+) + \pi_e\frac{\partial E}{\partial \bm{x}}
\end{align} \label{eq:corner_conditions}
where $\pi_e$ is a scalar Lagrange multiplier that can be calculated in terms of the states and costates. It should be noted that there is no discontinuous jump in the Hamiltonian due to the time regularization. Also, the primer vector, Eq. (\ref{eq:primer_vector}), and reduced switching function, Eq. (\ref{eq:switch}), are continuous across eclipse entry and exit \cite{wangtopputo}. Some recent works have successfully applied this traditional approach, also referred to as ``interior point constraints'', to low-thrust trajectory optimization with eclipsing \cite{wangtopputo,cerf2018,PONTANI2022549}. This traditional approach can be cumbersome, though, because it requires a trajectory propagation method that detects and stops exactly at every eclipse entry and exit. The eclipsing control constraint can also be incorporated into the dynamics with a smooth Heaviside approximation \cite{FerrierEpenoy}, which allows for straightforward trajectory propagation. The smooth, approximate Heaviside function approach will yield identical results to the discontinuous, true Heaviside function approach in the limit as smoothing parameters approach zero. The focus of this work is the averaged dynamics model, which does not exhibit discontinuous costate jumps, but does need to incorporate the averaged effect of these jumps into the averaged equations of motion. In the next section, the averaged dynamics are introduced, which use a novel application of the Leibniz integral rule to incorporate the effect of the costate jumps into the averaged equations of motion.

\section{Averaged Dynamics} \label{sec:averaging}

\subsection{Basic method} \label{sec:averaging_basic}

The averaging approach used in this work follows that of Edelbaum et al. \cite{edelbaum} and Kechichian \cite{kechichian}. The process starts with the unaveraged Hamiltonian including the embedded optimal control. Next, the average is performed over one period in time, and this averaged Hamiltonian is used to calculate the equations of motion for the averaged states and costates. It should be noted that the optimal control does not need to be re-derived after averaging the Hamiltonian, because there are no direct constraints on the control \cite{geffroy}. The averaging integral is converted from time to the fast moving true longitude, \textit{L}. All other variables are held constant within the averaging integral. This process removes the fast periodic changes in the state and costate, leaving only the secular motion. Gaussian quadrature \cite{quadrature} is used to approximate the integral and all other equations derived from it. The time averaged Hamiltonian, $\Tilde{H}$, over one period is defined here and rearranged slightly:
\begin{align}
\begin{split}
    \tilde{H} &= \frac{1}{T_0}\int_{-T_0/2}^{T_0/2} H dt = \frac{1}{T_0} \int_{-\pi}^{\pi} \frac{1}{\Dot{L}} H dL = \frac{1}{2 \pi} \int_{-\pi}^{\pi}  \frac{n}{\Dot{L}} H dL \\
    &\longrightarrow \tilde{H} = \frac{1}{2 \pi} \int_{-\pi}^{\pi}  s H dL, \ \ \ s = \frac{n}{\Dot{L}} 
\end{split}
\end{align}

\noindent where a new term, $s$, is defined for convenience and is only a function of the state variables, not the costates. The state and costate dynamics are calculated using the chain rule:
\begin{align}
\begin{split}
    \Dot{\text{\textbf{x}}} &= \frac{\partial}{\partial \bm{\uplambda}} \tilde{H} = \frac{\partial}{\partial \bm{\uplambda}} \left[\frac{1}{2 \pi} \int_{-\pi}^{\pi}  s  H  dL \right]= \frac{1}{2 \pi} \int_{-\pi}^{\pi}  s  \frac{\partial H}{\partial \bm{\uplambda}} dL 
\end{split}
\end{align}
\noindent \begin{align}
\begin{split}
    \Dot{\bm{\uplambda}} &= -\frac{\partial}{\partial \text{\textbf{x}}} \tilde{H} = -\frac{\partial}{\partial \text{\textbf{x}}} \left[\frac{1}{2 \pi} \int_{-\pi}^{\pi}  s  H  dL \right]= -\frac{1}{2 \pi} \int_{-\pi}^{\pi}  \frac{\partial}{\partial \text{\textbf{x}}}[s  H]  dL \\
    &\longrightarrow \bm{\Dot{\uplambda}} =  -\frac{1}{2 \pi} \int_{-\pi}^{\pi}  \left[\frac{\partial s}{\partial \text{\textbf{x}}} H + s \frac{\partial H}{\partial \text{\textbf{x}}}\right]  dL
\end{split}
\end{align}

\noindent noting that Leibniz terms are absent because the bounds on the integral are fixed. Mass is now included for notational convenience as a 9\textsuperscript{th} state in \textbf{x}:
\begin{align}
    \text{\textbf{x}} &= [p, f, g, h, k, L, t, \alpha, m]^{\intercal}, \ \ \bm{\uplambda} = [\lambda_p, \lambda_f, \lambda_g, \lambda_h, \lambda_k, \lambda_L, \lambda_t, \lambda_{\alpha}, \lambda_m]^{\intercal}
\end{align}

The \textit{n}-point Gaussian quadrature formula is used to approximate the integrals:
\begin{align}
    \int_a^b f(x) dx \approx \frac{b-a}{2}\sum_{i=1}^n w_i f\left(\frac{b-a}{2}z_i + \frac{b+a}{2}\right) \label{eq:quadrature}
\end{align}

\noindent where the weights, $w_i$, and node points, $z_i$, are calculated from Legendre polynomials \cite{quadrature}. The dynamics equations are then calculated using quadrature:
\begin{align}
\begin{split}
    \Dot{\text{\textbf{x}}} &= \frac{1}{2 \pi} \int_{-\pi}^{\pi}  s  \frac{\partial H}{\partial \bm{\uplambda}}  dL \approx \left. \frac{1}{2} \sum_{i=1}^n w_i \left[ s \frac{\partial H}{\partial \bm{\uplambda}} \right] \right\vert_{L_i}, \ \ \ \ \ L_i=\pi z_i \label{eq:x_dot_basic}
\end{split}
\end{align}
\noindent \begin{align}
\begin{split}
    \Dot{\bm{\uplambda}} &= -\frac{1}{2 \pi} \int_{-\pi}^{\pi}  \left[\frac{\partial s}{\partial \text{\textbf{x}}} H + s \frac{\partial H}{\partial \text{\textbf{x}}}\right]  dL \approx \left. -\frac{1}{2} \sum_{i=1}^n w_i \left[\frac{\partial s}{\partial \text{\textbf{x}}} H + s \frac{\partial H}{\partial \text{\textbf{x}}}\right] \right\vert_{L_i} \label{eq:lam_dot_basic}
\end{split}
\end{align}

The number of quadrature points required, \textit{n}, to make the approximation sufficiently accurate depends on the complexity of the unaveraged dynamics within the given period of averaging. If there are sections of the given period with and without thrusting, then \textit{n} must be large enough to accurately capture the proportion and relative locations of the two. This common \cite{geffroy,mazzini}, but brute force approach to calculating the averaged equations of motion can sometimes require such large \textit{n} that the computational savings of averaged propagation are significantly reduced. Alternatively, a multi-arc approach improves the accuracy of approximation without requiring excessively large \textit{n}.

\subsection{Multi-arc averaging} \label{sec:averaging_multi}

\begin{figure}[h!]
    \centering
    \includegraphics[scale=0.4]{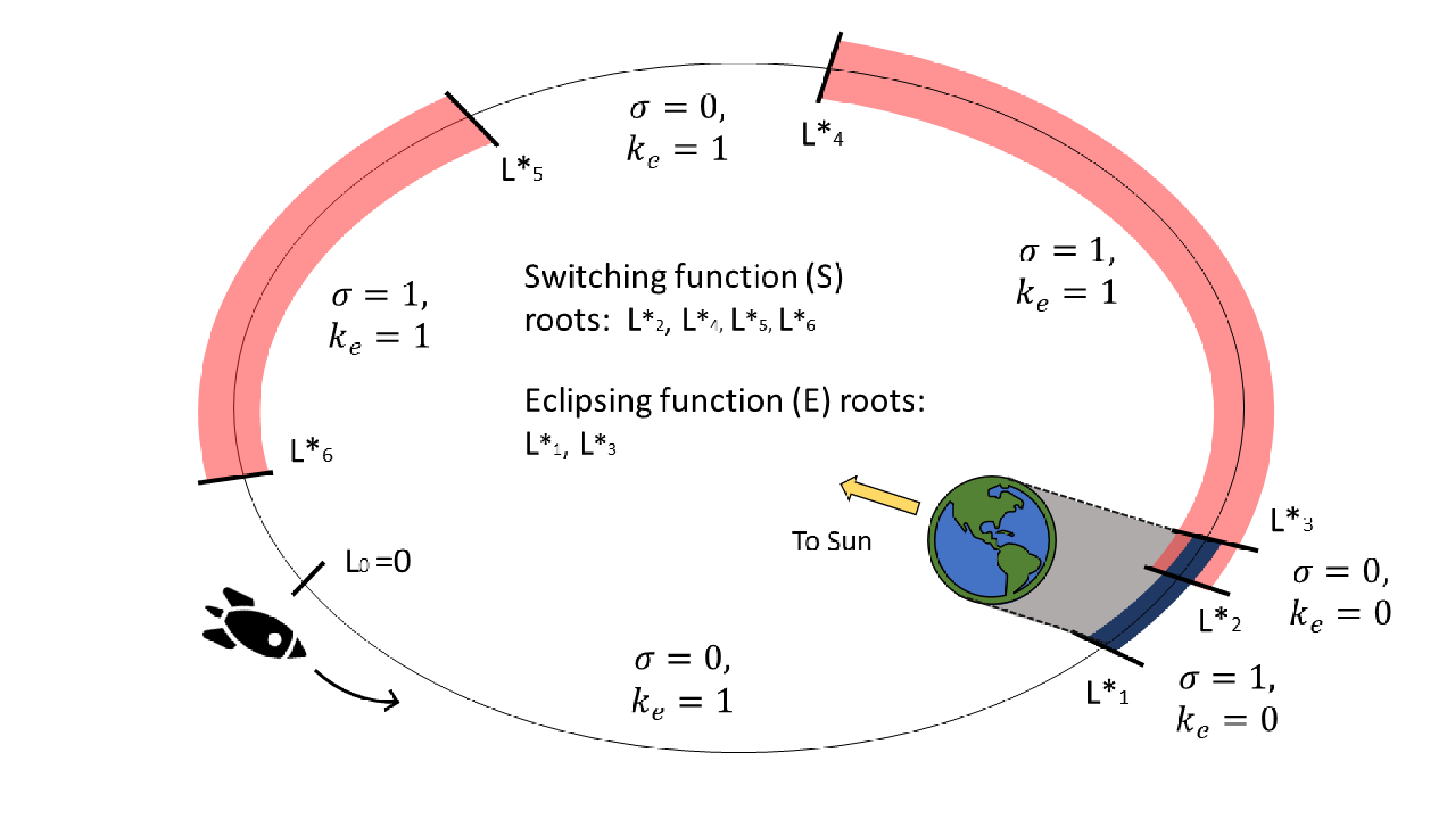}
    \caption{Diagram of multi-arc breakdown for a hypothetical spacecraft period.}
    \label{fig:multiarc}
\end{figure}

Introduced by Tarzi \cite{tarzi}, the multi-arc approach to bang-bang optimal control separates the averaging integral into different regions of thrust on/off depending on the roots of the eclipsing and switching function. Figure \ref{fig:multiarc} shows a diagram of a single period with 6 sections. Each section, or ``arc'', is treated separately in the averaging integral using the exact boundaries defined by the switching and eclipsing function roots. The segmenting of the averaging integral mitigates the need for a large quadrature number to accurately capture proportions of thrusting and coasting. The arcs with no thrusting due to shadow are treated separately from the arcs that have no thrusting due to the bang-bang optimal control law. Each arc uses fixed values of either 0 or 1 for $k_e$ and $\sigma$, depending on the sign of the eclipsing and switching functions, respectively. The multi-arc averaged dynamics approximate the true bang-bang control of the unaveraged minimum fuel problem.

\begin{figure}[h!]
\centering
\begin{minipage}{1\textwidth}
  \includegraphics[width=1\linewidth]{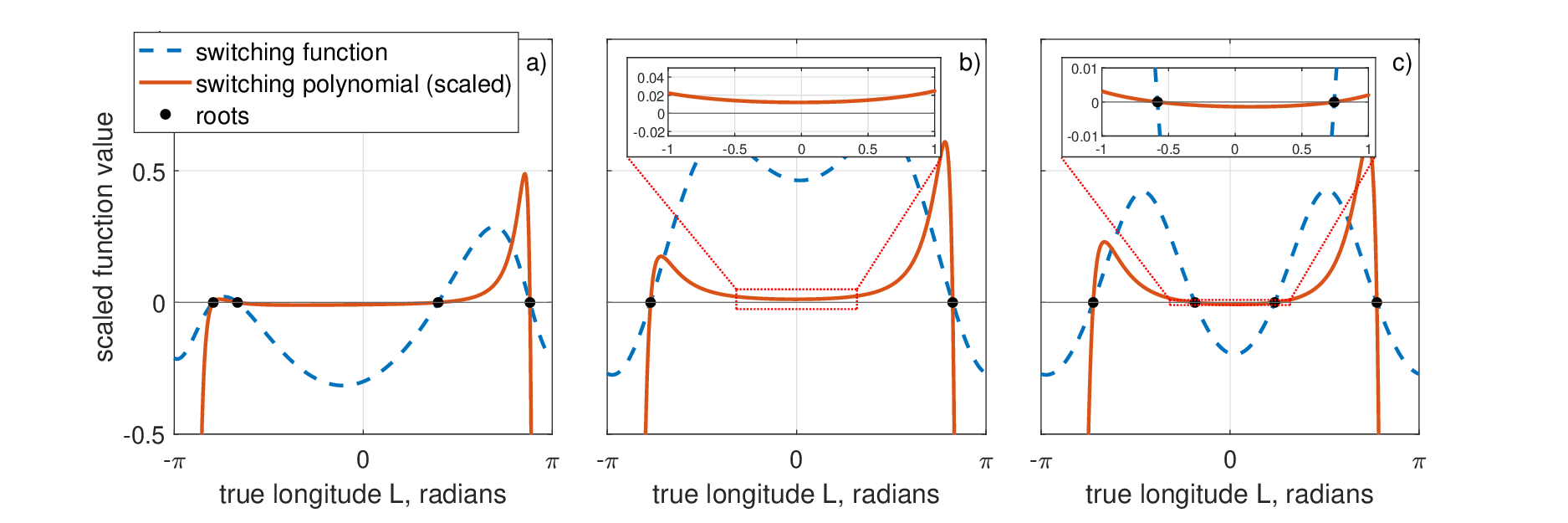}
\end{minipage}
\caption{Switching function and switching polynomial over averaging period at (a: left) beginning of transfer, (b: middle) middle of transfer, and (c: right) end of transfer.}
\label{fig:switch_poly}
\end{figure}

The \textit{M} roots of the switching function from Eq. (\ref{eq:switch}) are computed exactly using a new approach. When the switching function is fully expanded in terms of \textbf{B}, a sixth-order polynomial of tan$(L/2)$ can be shown to have the same roots when all variables except $L$ are held constant within the averaging integral. The switching function roots are calculated at every integration step when using averaged dynamics, so the use of a simpler ``switching polynomial'', $S_p$, speeds up the overall process. A symbolic manipulator is used to show equivalency between the roots of the switching function in Eq. (\ref{eq:switch}) and those of $S_p$ in Eq. (\ref{eq:poly}):
\begin{align}
\begin{split}
    S_p &=  \beta_1 \tan^6(L/2) + \beta_2 \tan^5(L/2) + \beta_3 \tan^4 (L/2) + \beta_4 \tan^3 (L/2) + \beta_5 \tan^2 (L/2) + \beta_6 \tan(L/2) + \beta_7 \label{eq:poly}
\end{split}
\end{align}

\noindent where all non-longitude variables are combined and shown as $\bf{\beta}$ functions here for brevity. Full expressions $\bf{\beta}$ are shown in the appendix. This polynomial has the same roots as the original switching function, but is simpler for root-finding. An example of the switching function and switching polynomial over an averaging period is shown in Fig. \ref{fig:switch_poly} for beginning, middle, and end of the first transfer shown in Section \ref{sec:results}.

The form of  Eq. (\ref{eq:poly}) implies that a maximum of 6 relevant roots exist, depending on the $\beta$ coefficients, noting that any complex roots are not relevant. Therefore, a maximum number of three separate arcs of thrusting are possible within a single period in the averaging context for a minimum fuel problem. This theoretical result, consistent with the numerical results of \cite{taheri2020many}, implies a maximum number of optimal thrusts per revolution in the low-thrust limit. Previous works investigating the structure of the switching function have set limits of 8 unique roots per period specifically during coasting in unaveraged minimum fuel problems \cite{Pan2012,Jamison2010} with lower limits when additional constraints are included. The new averaged dynamics-based upper limit of 6 roots provides further insight on Edelbaum's question of ``how many impulses'' are optimal for transfers between orbits when averaging assumptions hold \cite{edelbaum_howmanyimpulses,saloglu}. Even in direct optimization, the number of impulses or thrust arcs per revolution has a large impact on run time and resolution of the problem. A theoretical upper bound on the number of thrusts per revolution is therefore a significant result for practical mission design. This theoretical result was independently developed and first presented to the public by the authors in January, 2024 \cite{lifset}. A similar result was presumably developed simultaneously, and subsequently published by Wu et al. in March 2024 \cite{wu_avg}.

The eclipsing function is much simpler than the switching function, with at most two roots for a single revolution. The roots are calculated using a Newton-Raphson routine with initial guesses found analytically from a cylindrical shadow model \cite{vallado2001fundamentals}. After the roots of the switching and eclipsing functions are calculated, the averaged Hamiltonian is broken down into a sum across $M$ arcs:
\begin{align}
    \tilde{H} &= \frac{1}{2 \pi} \int_{-\pi}^{ \pi}  s H dL = \frac{1}{2 \pi} \sum_{j=1}^M\int_{L^*_{j-1}}^{L^*_j}  s_j H_j dL \label{eq:H_int}
\end{align}

\noindent where the $L^*$ variables specifically refer to roots of the switching and eclipsing functions, which are functions of the other state variables.

The dynamics derivative calculations  require the Leibniz integral rule, because the integral limits are now functions of the differentiation variables. In general form, the Leibniz integral rule \cite{Protter1985} is
\begin{align}
    \frac{d}{d x} \left[\int_{a(x)}^{b(x)} f(x,t) dt \right] &= f(x,b(x)) \frac{d}{d x} b(x) - f(x,a(x)) \frac{d}{d x} a(x) + \int_{a(x)}^{b(x)} \frac{d f(x,t)}{d x} dt
\end{align}

\noindent This rule is used on Eq. (\ref{eq:H_int}) to calculate
\begin{align}
\begin{split}
    \Dot{\bm{\uplambda}} &= -\frac{\partial \Tilde{H}}{\partial \text{\textbf{x}}} = -\frac{1}{2 \pi} \sum_{j=1}^M \frac{\partial}{\partial \text{\textbf{x}}} \int_{L^*_{j-1}(\text{\textbf{x}},\bm{\uplambda})}^{L^*_j(\text{\textbf{x}},\bm{\uplambda})}  s_j H_j  dL = - \frac{1}{2 \pi} \sum_{j=1}^M \left\{ \left. \int_{L^*_{j-1}(\text{\textbf{x}},\bm{\uplambda})}^{L^*_j(\text{\textbf{x}},\bm{\uplambda})}  \frac{\partial}{\partial \text{\textbf{x}}} [s_j H_j]  dL + \frac{d L^*}{d \text{\textbf{x}}}s_j H_j \right\vert^{L^*_j}_{L^*_{j-1}} \right\} \\
    &= - \frac{1}{2 \pi} \sum_{j=1}^M \left\{ \left. \frac{L^*_j - L^*_{j-1}}{2} \sum^{n_j}_{i=1} w_i \left. \left[ \frac{\partial s_j}{\partial \text{\textbf{x}}}H_j + s_j\frac{\partial H_j}{\partial \text{\textbf{x}}}\right] \right|_{L_i} + \frac{d L^*}{d \text{\textbf{x}}}s_j H_j \right\vert^{L^*_j}_{L^*_{j-1}} \right\},  \ \ \ L_i = \left[\frac{L^*_j - L^*_{j-1}}{2}\right]z_i + \frac{L^*_j + L^*_{j-1}}{2} \label{eq:lam_dot}
\end{split}
\end{align}
\noindent \begin{align}
\begin{split}
    \Dot{\text{\textbf{x}}} &= \frac{\partial \Tilde{H}}{\partial \bm{\uplambda}} = \frac{1}{2 \pi} \sum_{j=1}^M \frac{\partial}{\partial \bm{\uplambda}} \int_{L^*_{j-1}(\text{\textbf{x}},\bm{\uplambda})}^{L^*_j(\text{\textbf{x}},\bm{\uplambda})}  s_j H_j  dL = \frac{1}{2 \pi} \sum_{j=1}^M \left\{ \left. \int_{L^*_{j-1}(\text{\textbf{x}},\bm{\uplambda})}^{L^*_j(\text{\textbf{x}},\bm{\uplambda})}  \frac{\partial}{\partial \bm{\uplambda}} [s_j H_j]  dL + \frac{d L^*}{d \bm{\uplambda}}s_j H_j \right\vert^{L^*_j}_{L^*_{j-1}} \right\} \\
    &= \frac{1}{2 \pi} \sum_{j=1}^M \left\{ \left. \frac{L^*_j - L^*_{j-1}}{2} \sum^{n_j}_{i=1} w_i \left. \left[ s_j\frac{\partial H_j}{\partial \bm{\uplambda}}\right] \right|_{L_i} + \frac{d L^*}{d \bm{\uplambda}}s_j H_j \right\vert^{L^*_j}_{L^*_{j-1}} \right\},  \ \ \ L_i = \left[\frac{L^*_j - L^*_{j-1}}{2}\right]z_i + \frac{L^*_j + L^*_{j-1}}{2} \label{eq:x_dot}
\end{split}
\end{align}

\noindent The number of quadrature points required is reduced compared to the brute-force averaging approach. In this multi-arc model, the number of quadrature points used per arc is defined to be a function of the length of each arc:
\begin{align}
    n_j &= q \ [1 + 2 \ \text{round}(L^*_j - L^*_{j-1})] \label{eq:quadrature_number}
\end{align}

\noindent where $L^*$ takes units of radians and \textit{q} is an integer-valued tuning parameter. 

\begin{figure}[h!]
    \centering
    \includegraphics[scale=0.4]{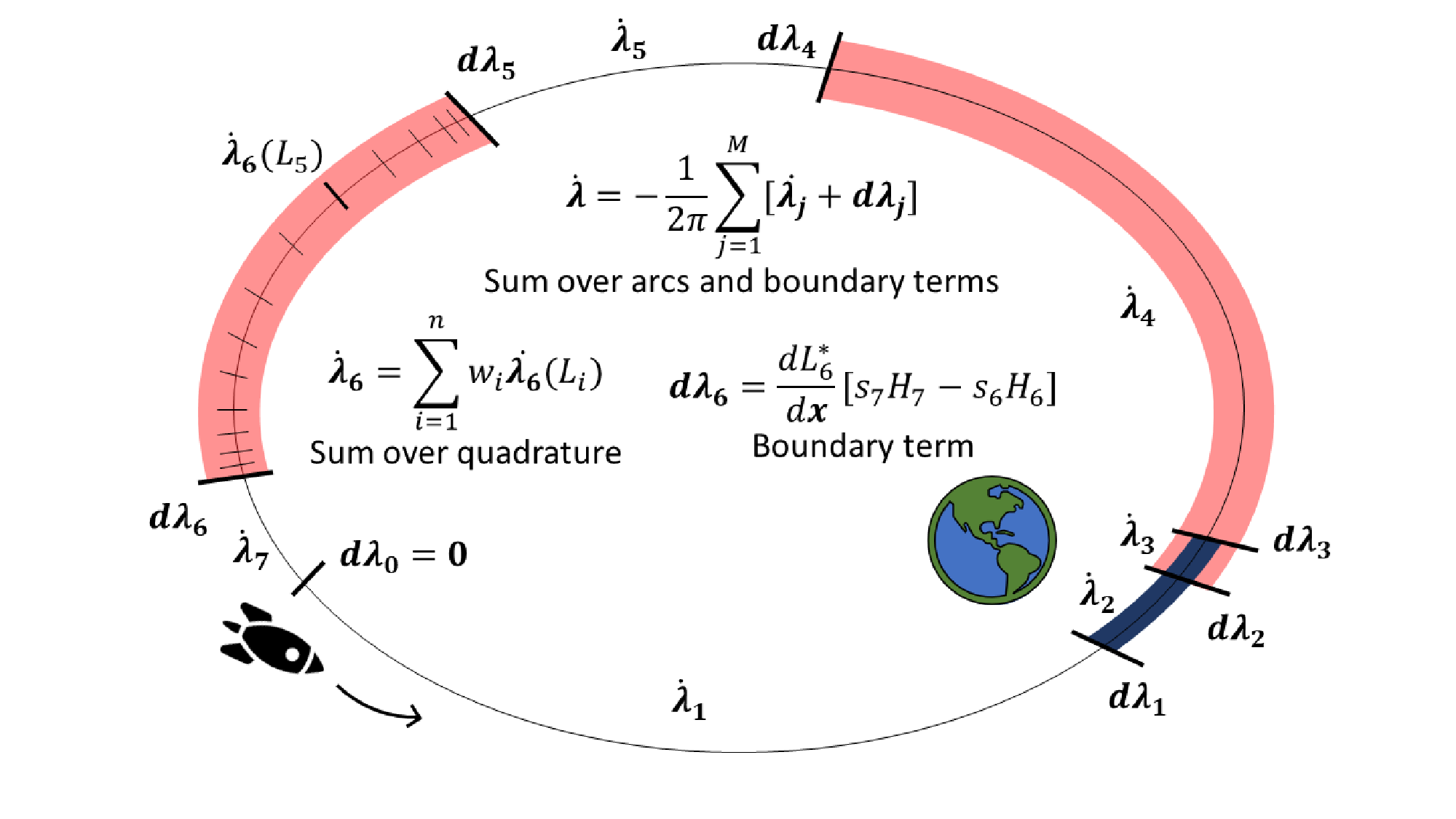}
    \caption{Contributions to averaged costate dynamics in a hypothetical multi-arc period. The equation shown is simplified with separated integral and Leibniz boundary terms. }
    \label{fig:multiarc_math}
\end{figure}

The Leibniz integral rule terms in Eqs. (\ref{eq:lam_dot}) and (\ref{eq:x_dot}) capture the effect of costate jumps due to the eclipse constraint in the dynamics, and naturally smooth them out over the averaging period. These terms are functionally equivalent to the Dirac delta terms used by Mazzini\cite{mazzini}. Note also that the Leibniz terms are related to the additional terms attributed to sensitivities propagated across bang-bang switches (see Eq. (29) in Ref. \cite{russell2007primer}). Figure \ref{fig:multiarc_math} shows how the quadrature points and Leibniz boundary terms contribute to the costate dynamics for the same hypothetical multi-arc period. For calculation of these terms, both the derivative and the $s_j$ or $H_j$ are evaluated at true longitude roots of the eclipsing and switching functions. The partial derivative of the switching boundary terms are computed assuming the switching function constraint is satisfied:
\begin{align}
\begin{split}
    S &= S(\text{\textbf{x}}, \bm{\uplambda}, L^*_j(\text{\textbf{x}}, \bm{\uplambda})) = 0 \longrightarrow \frac{d S}{d \text{\textbf{x}}} = \frac{\partial S}{\partial \text{\textbf{x}}} + \frac{\partial S}{\partial L} \frac{d L^*_j}{d \text{\textbf{x}}} = 0 \\
    &\longrightarrow \frac{d L^*_j}{d \text{\textbf{x}}} = - \frac{\partial S}{\partial \text{\textbf{x}}}/\frac{\partial S}{\partial L}, \ \ \ \frac{d L^*_j}{d \bm{\uplambda}} = - \frac{\partial S}{\partial \bm{\uplambda}}/\frac{\partial S}{\partial L} \label{eq:DSDL}
\end{split}
\end{align}

\noindent For the eclipse entrance/exit points, the situation is simpler, because the eclipsing function is not a function of the costate:
\begin{align}
\begin{split}
    E &= E(\text{\textbf{x}}, L^*_j(\text{\textbf{x}})) = 0 \longrightarrow \frac{d E}{d \text{\textbf{x}}} = \frac{\partial E}{\partial \text{\textbf{x}}} + \frac{\partial E}{\partial L} \frac{d L^*_j}{d \text{\textbf{x}}} = 0 \\
    &\longrightarrow \frac{d L^*_j}{d \text{\textbf{x}}} = - \frac{\partial E}{\partial \text{\textbf{x}}}/\frac{\partial E}{\partial L}, \ \ \ \frac{d L^*_j}{d \bm{\uplambda}} = 0 \label{eq:DEDL}
\end{split}
\end{align}

\noindent The switching, Eq. (\ref{eq:DSDL}), or eclipsing, Eq. (\ref{eq:DEDL}), version is used depending on whether the relevant integral boundary is due to a switch or an eclipse.


\subsection{Eclipsing singularity} \label{sec:singularity}

The new Leibniz boundary terms in Eq. (\ref{eq:lam_dot}) smoothly incorporate the costate jumps into the averaged equations of motion. As eclipsing arcs get smaller and approach zero length, though, these extra terms increase in magnitude. Unfortunately, they become singular in the limit when the arc disappears. The averaged costate jumps used in Mazzini's averaging model \cite{mazzini2} similarly become singular. The singularity can be shown mathematically by breaking down a hypothetical averaging period with an eclipse and no coasting arcs:
\begin{align}
\begin{split}
    \Dot{\bm{\uplambda}} &= -\frac{1}{2 \pi} \left\{ \int_{-\pi}^{L^*_{\text{in}}}  \frac{\partial}{\partial \text{\textbf{x}}} (s_1 H_1)  dL + \int_{L^*_{\text{in}}}^{L^*_{\text{out}}}  \frac{\partial}{\partial \text{\textbf{x}}} (s_2 H_2)  dL + \int_{L^*_{\text{out}}}^{\pi}  \frac{\partial}{\partial \text{\textbf{x}}} (s_1 H_1)  dL + \underbrace{\left[ \frac{d L^*_{\text{in}}}{d \text{\textbf{x}}} - \frac{d L^*_{\text{out}}}{d \text{\textbf{x}}} \right]}_{\text{term }1}\underbrace{[s_1 H_1 - s_2 H_2]}_{\text{term }2} \right\} \label{eq:singular_lam_dot}
\end{split}
\end{align}

\noindent where $(L^*_{\text{in}},L^*_{\text{out}})$ are the true longitude values for the entrance and exit of the eclipse, respectively. As an eclipsing arc approaches zero length, term 1 of Eq. (\ref{eq:singular_lam_dot}) approaches infinity, while term 2 stays finite. A trajectory propagated with unaveraged dynamics will only be affected by this singularity if it happens to cross a point where the eclipse arc is very small (approximately less than $0.01$ radians) and incorporates the eclipsing constraint in a discontinuous way \cite{wangtopputo}. However, a trajectory propagated with averaged dynamics will encounter this singularity if any part of the averaging integral contains a very small eclipsing arc, which occurs nearly every time eclipsing is present. A trajectory propagated with averaged dynamics and a variable step integrator will significantly slow down when approaching the singularity, and will not be able to propagate through it without enforcing a minimum step size. Additionally, a trajectory that is forced to propagate through the singularity will accumulate errors in the equations of motion. 

\begin{figure}[h!]
    \centering
    \includegraphics[scale=0.65]{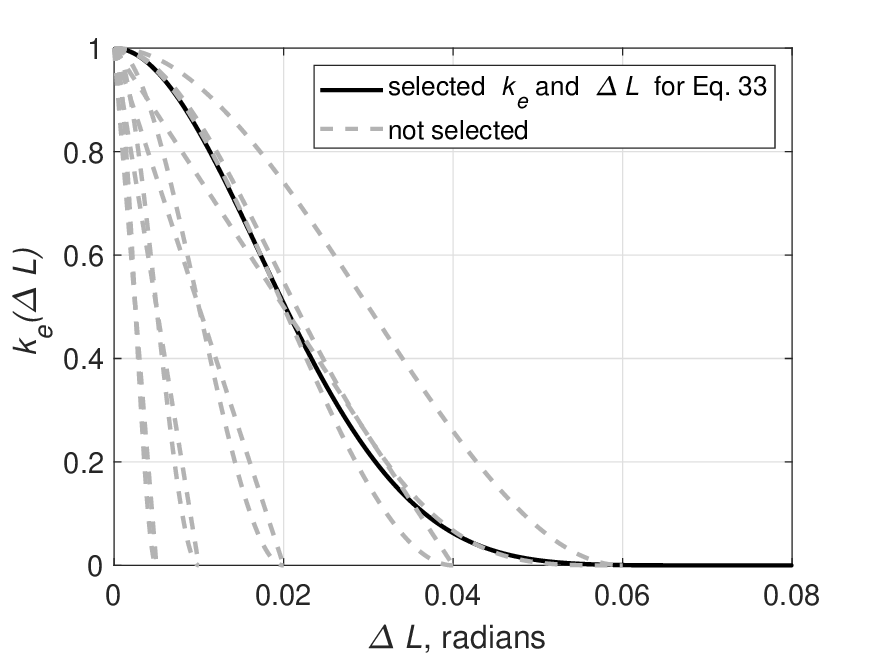}
    \caption{Updated minimum thrust value inside eclipse.}
    \label{fig:ke_fix}
\end{figure}

Note that no singularity is present in the physics of a realistic eclipse.  The true penumbra and the finite radius of a spacecraft body lead to a continuous transition between full and zero sunlight. Therefore, the singularity described here is an artifact of a convenient eclipsing model that is both approximate and discontinuous. The modeling singularity is removed through a novel and minor redefinition of the eclipsing constraint. As the eclipsing arc approaches zero length, the minimum (normalized) thrust inside the eclipse is slowly raised from $0$ to $1$:
\begin{align}
\begin{split}
    k_e &= 
    \begin{cases}
        1 & \text{if  } E < 0 \\
        0 & \text{if  } E > 0 \text{ and } \Delta L> \Delta L^\star\\
        k_e(\Delta L) & \text{if  } E > 0 \text{ and } \Delta L< \Delta L^\star
    \end{cases}
\end{split}
\end{align}

\noindent where $\Delta L$ is the eclipsing arc size in radians. This redefinition causes term 2 of Eq. (\ref{eq:singular_lam_dot}) to approach zero as the eclipsing arcs become small, effectively removing the singularity by canceling with the increase of term 1. A small $\Delta L^\star$ is preferred, so that the $k_e(\Delta L)$ redefinition has a minimal effect on the propagated trajectory. However, a larger $\Delta L^\star$ with a gradually decreasing $k_e(\Delta L)$ function results in a faster propagation. A variable step propagator will always slow down as eclipsing arcs disappear, but the extent of slow-down is affected by this choice. Trade studies are conducted on the function $k_e(\Delta L)$ and the arc size where the fix begins, $\Delta L^\star$ to balance these two effects. Linear through 12th order polynomials are tested with fixed points at $k_e(0) = 1$ and $k_e(\Delta L^\star) = 0$. The following value of $\Delta L^\star$ and polynomial for $k_e$ are selected to balance compute time and realism:
\begin{align}
    \Delta L^\star &= 0.08 \text{ radians} \\
    k_e(\Delta L) &= \frac{1}{256}[15625 \Delta L ^3-1875 \Delta L ^2+4]^4 \label{eq:ke_fix}
\end{align}
The function is shown in Fig. \ref{fig:ke_fix} along with other tested polynomials. 

Wang and Topputo \cite{wangtopputo} use a similar technique of varying the floor of a discontinuous eclipsing constraint in their approach to the indirect method. Their model is unaveraged, but converges better when the eclipsing constraint floor is slowly lowered using a continuation scheme. The fix shown here is motivated primarily by a singularity in the equations of motion, but improves convergence as well. The overall flow of calculating the multi-arc averaged dynamics is shown in Alg. \ref{alg:cap}. 

\begin{algorithm}[h!]
\footnotesize
\caption{Multi-arc Averaged Dynamics}\label{alg:cap}
\begin{algorithmic}
\State \textbf{Inputs:} $\tau$, \textbf{x}, $\bm{\uplambda}$ \\
\State Solve for exact roots of switching and eclipsing functions
\State \hskip1.5em Calculate $L^*$ such that $S_p(L^*) = 0$
\State \hskip1.5em Calculate $L^*$ such that $E(L^*) = 0$
\State \hskip1.5em Discard imaginary roots
\If{No exact roots found}
    \State Calculate $\Dot{\text{\textbf{x}}}$ and $\bm{\Dot{\uplambda}}$ with single arc averaging equations \ref{eq:x_dot_basic} and \ref{eq:lam_dot_basic} \\
    \hskip1.5em \textbf{return}
\EndIf
\State Sort exact roots, $L_i^*$, to determine the \textit{M} arcs
\State \textbf{for} arc \textit{j} = 1 to \textit{M}
\State \hskip1.5em \textbf{if} S < 0
\State \hskip1.5em \hskip1.5em $\sigma = 1$
\State \hskip1.5em \textbf{else}
\State \hskip1.5em \hskip1.5em $\sigma = 0$
\State \hskip1.5em \textbf{end if}
\State \hskip1.5em \textbf{if} E < 0
\State \hskip1.5em \hskip1.5em $k_e = 1$
\State \hskip1.5em \textbf{else}
\State \hskip1.5em \hskip1.5em \textbf{if} $\Delta L > \Delta L^\star$
\State \hskip1.5em \hskip1.5em \hskip1.5em $k_e = 0$
\State \hskip1.5em \hskip1.5em \textbf{else}
\State \hskip1.5em \hskip1.5em \hskip1.5em $k_e = k_e(\Delta L)$ from Eq. (\ref{eq:ke_fix})
\State \hskip1.5em \hskip1.5em \textbf{end if}
\State \hskip1.5em \textbf{end if}
\State \hskip1.5em Calculate $n_j$ with Eq. (\ref{eq:quadrature_number})
\State \hskip1.5em Calculate $\Dot{\text{\textbf{x}}}_j$ and $\bm{\Dot{\uplambda}}_j$ with quadrature and boundary terms \Comment{The bracketed term in Eqs. (\ref{eq:lam_dot}) and (\ref{eq:x_dot}})
\State \textbf{end for}
\State Calculate $\Dot{\text{\textbf{x}}}$ and $\bm{\Dot{\uplambda}}$ by summing $\Dot{\text{\textbf{x}}}_j$ and $\bm{\Dot{\uplambda}}_j$ over \textit{M} arcs \\
\textbf{return} \\ \\ 
\State \textbf{Outputs:} $\Dot{\text{\textbf{x}}}$, $\bm{\Dot{\uplambda}}$
\end{algorithmic}
\end{algorithm}

\section{Variational Equations} \label{sec:variational}

The state transition matrix (STM) of the full augmented state vector is useful for a variety of applications, including solving boundary value problems. Previous works using averaged dynamics either do not mention the STM when solving the BVP \cite{geffroy,mazzini} or use a finite difference method to calculate it \cite{tarzi}. In this work, the STM, $\Phi$, is integrated along with the combined states and costates using the variational equations. While tedious to derive, this approach is faster and more accurate. The variational equations for averaged dynamics are explicitly provided here for the first time in the literature. It is possible that prior authors have utilized the variational equations to compute STMs and solve BVPs, but details of their calculation, including methods to avoid singularities, have not yet appeared in the averaging literature. The STM is integrated using the typical approach with the variational equations:
\begin{align}
\begin{split}
    \Phi(\tau) &\equiv \frac{d \text{\textbf{y}}(\tau)}{d \text{\textbf{y}}_0} \\
    \frac{d}{d \tau} \Phi(\tau) &= \frac{d \Dot{\text{\textbf{y}}}(\tau) }{d \text{\textbf{y}}(\tau)} \Phi(\tau), \ \ \ \Phi_0 = \textbf{I} _{18\times18}
\end{split}
\end{align}

\noindent where \textbf{y} is the augmented state and costate vector, and a flipped vector is defined for notational convenience:
\begin{align}
\begin{split}
    \text{\textbf{y}} &= [\text{\textbf{x}}^{\intercal}, \bm{\uplambda}^{\intercal}]^{\intercal}, \ \ \ \text{\textbf{y}}^{\text{F}} = [\bm{\uplambda}^{\intercal}, \text{\textbf{x}}^{\intercal}]^{\intercal}
\end{split}
\end{align}

The averaged dynamics require extra attention to account for the Leibniz integral terms. The roots of the eclipsing and switching functions, $L^*$ are functions of the state variables, and must be considered in the total derivative. Therefore, the Jacobian of the dynamics is calculated using the quadrature formulation of Eqs. (\ref{eq:x_dot}) and (\ref{eq:lam_dot}):
\begin{align}
\begin{split}
    \frac{d \Dot{\Tilde{\text{\textbf{y}}}}(\tau) }{d \Tilde{\text{\textbf{y}}}(\tau)} &= \frac{1}{2 \pi} \sum_{j=1}^M \left\{ \frac{1}{2} \left[\frac{\partial L^*_j}{\partial \Tilde{\text{\textbf{y}}}} - \frac{\partial L^*_{j-1}}{\partial \Tilde{\text{\textbf{y}}}}\right] \sum_{i=1}^{n_j} w_i \left.\frac{\partial }{\partial \Tilde{\text{\textbf{y}}}^F} (s_j H_j)\right|_{L_i} \right.\\
    & \left.+ \frac{1}{2} \left[L^*_j - L^*_{j-1}\right] \sum_{i=1}^{n_j} w_i \left. \left[ \frac{\partial^2 }{\partial \Tilde{\text{\textbf{y}}} \partial \Tilde{\text{\textbf{y}}}^F} (s_j H_j)\right|_{L_i} + \frac{d L_i}{d \Tilde{\text{\textbf{y}}}} \left.\frac{\partial^2 }{\partial \Tilde{L} \partial \Tilde{\text{\textbf{y}}}^F} (s_j H_j)\right|_{L_i}\right]  \right. \\
    &+ \left. \left. \left[ \frac{\partial^2 L^*}{\partial \Tilde{\text{\textbf{y}}} \partial \Tilde{\text{\textbf{y}}}^F} + \frac{d L^*}{d \Tilde{\text{\textbf{y}}}}\frac{\partial^2 L^*}{\partial \Tilde{L} \partial \Tilde{\text{\textbf{y}}}^F}\right]s_j H_j \right|_{L^*_{j-1}}^{L^*_j} + \left. \frac{d L^*}{d \Tilde{\text{\textbf{y}}}^F}\left[ \frac{\partial}{\partial \Tilde{\text{\textbf{y}}}}(s_j H_j) + \frac{d L^*}{\partial \Tilde{\text{\textbf{y}}}}\frac{\partial}{\partial \Tilde{L}} (s_j H_j) \right] \right|_{L^*_{j-1}}^{L^*_j} \right\} \\
    \frac{d L_i}{d \Tilde{\text{\textbf{y}}}} &= \frac{1}{2} \left[\frac{d L^*_j}{d \Tilde{\text{\textbf{y}}}} - \frac{d L^*_{j-1}}{d \Tilde{\text{\textbf{y}}}}\right]z_i + \frac{1}{2} \left[\frac{d L^*_j}{d \Tilde{\text{\textbf{y}}}} + \frac{d L^*_{j-1}}{d \Tilde{\text{\textbf{y}}}}\right] \label{eq:variational}
\end{split}
\end{align}

\noindent The second derivatives of $L^*$ are calculated by taking derivatives of Eqs. (\ref{eq:DEDL}) and (\ref{eq:DSDL}). These variational equations and the averaged dynamics in general will be validated in the next section with an optimal transfer example.

\section{Example Optimal Trajectories} \label{sec:results}

Before the results are shown, an approximate unaveraged dynamics model is introduced. Discontinuous bang-bang control and eclipsing constraints have been used with unaveraged optimal control by Wang and Topputo \cite{wangtopputo}, Russell \cite{russell2007primer}, and others. This approach typically requires a propagator with integrator events, complicated STM calculations, and explicit costate jump calculations. A common approach to approximate the unaveraged  model uses smooth Heaviside approximations for bang-bang control and the eclipsing constraint, and eliminates the need for event finding or jumps in the costate STMs. This approximate model is introduced here solely for validation of the averaged model. The approximate Heavisides are set to be extremely steep and as close to true bang-bang control as possible to make the comparison to the averaged dynamics sufficiently accurate. The use of approximate Heavisides for smooth bang-bang control and eclipsing in averaged dynamics is beyond the current scope. For more details on approximate Heaviside implementations in optimal control, see \cite{taheri} and \cite{bertrandepenoy}.

Two examples are presented: A GTO to GEO transfer is shown to validate the averaged model through comparison with unaveraged dynamics. The averaged and unaveraged transfers are optimized separately to show how the optimal averaged trajectory can be used as a proxy for the optimal unaveraged trajectory. A GTO to GEO optimal transfer is then shown using only averaged dynamics. This latter example is almost 500 revolutions and is intractable to practically optimize with unaveraged dynamics.


\subsection{Approximate Unaveraged Dynamics} \label{sec:approx_unaveraged}

The approximate unaveraged dynamics make trajectory propagation faster and smoother by substituting smooth Heaviside approximations for bang-bang control and the eclipsing constraint \cite{taheri}. The optimal control derivation and resulting equations of motion are very similar to those derived in Section \ref{sec:optimal_control} with only a few adjustments. 

First, the eclipsing constraint, Eq. (\ref{eq:ke}), is changed to a smooth Heaviside approximation:

\begin{align}
    k_e = \frac{1}{2}\left(1 - \frac{E}{\sqrt{E^2+\epsilon_{E}^2}}\right) \label{eq:ke_smooth}
\end{align}
where $\epsilon_{E}$ is a smoothing parameter. As $\epsilon_{E}$ is reduced to zero, Eq. (\ref{eq:ke_smooth}) approximates the discontinuous eclipsing constraint, Eq. (\ref{eq:ke}). It should be noted that the smooth Heaviside approximations asymptotically approach zero and one, but the averaged model sets the thrust to be equal to exactly zero or one. As a result of the smoothed eclipsing constraint, the costate jumps are also smoothed, but will approach discrete jumps as $\epsilon_{E}$ approaches zero. This L2-norm based Heaviside approximation is a an efficient alternative to the trigonometric form \cite{taheri_L2norm,Mall_L2norm,Heidrich_L2norm}. 

The cost function is changed to:
\begin{align}
    J &= \int_0^1 \left\{\alpha \frac{T_{\text{min}}}{c} + \alpha \frac{T_{\text{max}}-T_{\text{min}}}{c}k_e \left[\sigma - \epsilon_S\sqrt{\sigma - \sigma^2}\right] \right\}d\tau \label{eq:cost_smooth}
\end{align}
where $\epsilon_S$ is a smoothing parameter that is gradually reduced towards zero to approximate the true minimum fuel problem, Eq. (\ref{eq:PI}). The resulting Hamiltonian is:
\begin{align}
    H &= \alpha \frac{T_{\text{min}}}{c} + \alpha \frac{T_{\text{max}}-T_{\text{min}}}{c}k_e \left[\sigma - \epsilon_S\sqrt{\sigma - \sigma^2}\right] \nonumber \\
    &+ \alpha \bm{\uplambda}^{T} \text{\textbf{A}} + \alpha \bm{\uplambda}^{T} \text{\textbf{B}} \left[\hat{\bm{u}}  \frac{T}{m} +\gamma\right] - \alpha \lambda_m \frac{T}{c} \label{eq:hamiltonian_smooth}
\end{align}
By applying the smoothing approximation in the cost function definition, the derived optimal control law results in the same L2-norm based form of the Heaviside approximation used in Eq. (\ref{eq:ke_smooth}):
\begin{align}
    \sigma = \frac{1}{2}\left[1-\frac{S}{\sqrt{S^2 + \epsilon_S^2}} \right] \label{eq:sigma_smooth}
\end{align}
Besides the new definitions in Eq. (\ref{eq:ke_smooth}) and Eq. (\ref{eq:sigma_smooth}), all other equations from Section \ref{sec:optimal_control} are the same. When compared with averaged dynamics with bang-bang control and eclipsing shown in Sec. \ref{sec:averaging}, the smoothing parameters $[\epsilon_E,\epsilon_S]$ are set to the order of $10^{-5}$. The smooth Heaviside approximations introduced here are unrelated in mathematical structure with the Heaviside redefinition introduced in Sec. \ref{sec:singularity}. The two Heaviside adjustments are related in motivation, though, because they both benefit the simulation of trajectories.


\subsection{48-Revolution GTO to GEO example} \label{sec:GTOtoGEO}

A 48-revolution, optimal minimum-fuel GTO to GEO transfer is presented to compare the averaged and unaveraged dynamics models. The states and costates of the optimal averaged trajectory are directly compared with the optimal unaveraged trajectory for the first time in the literature. The transfer parameters and details are given in Tables \ref{tab:GTOtoGEO_BVP} and \ref{tab:GTOtoGEO_details}, respectively. Note that the parameters of this first example require an unrealistically high power level. However, they are chosen to simply validate the averaged dynamics model. A more practical transfer is shown in the later example. The Sun position and velocity data are calculated using the SPICE for MATLAB toolkit \cite{spice}, with a launch date of Apr 1, 2008 (260280065.0 ET). The trajectory is propagated using a Runge-Kutta 8(9) variable step integrator with relative and absolute tolerances set to $10^{-14}$. The averaged trajectory is propagated using $q=6$ in Eq. (\ref{eq:quadrature_number}). The unaveraged trajectory is propagated with smoothing parameters of $[\epsilon_E,\epsilon_S] = [3 \times 10^{-5}, 1 \times 10^{-5}]$ to approximate bang-bang control and compensate for scaling differences between the switching and eclipsing function amplitudes during the transfer. These values are selected by propagating an initial guess and iterating on the smoothing parameters until the steepness of $\sigma$ and $k_e$ look consistent. The full $J_2$ perturbation is included using modified equinoctial elements in the RTN frame \cite{kechichian}:
\begin{align}
\begin{split}
    \gamma_R &= -\frac{3 \mu J_2 R_{\text{E}}^2}{2 r^4} \left\{1 - \frac{12[h \sin(L) - k \cos(L)]}{[h^2 + k^2 + 1]^2}\right\} \\
    \gamma_T &= -\frac{12 \mu J_2 R_{\text{E}}^2}{r^4} \frac{[h \sin(L) - k \cos(L)][k \cos(L) + h \sin(L)]}{[h^2 + k^2 + 1]^2} \\
    \gamma_N &= -\frac{6 \mu J_2 R_{\text{E}}^2}{r^4} \frac{[k \sin(L) - h \cos(L)][-h^2 - k^2 + 1]}{[h^2 + k^2 + 1]^2} 
\end{split}
\end{align}

\noindent The boundary value problem is solved using MATLAB's FMINUNC to minimize the cost function, $Z$, defined as:
\begin{align}
    Z(\bm{\uplambda}_0) = ||\text{\textbf{x}}_f - \text{\textbf{x}}_{\text{GEO}}||^2 +\left[\lambda_{m,f}\right]^2+\left[\lambda_{L,f}\right]^2 \label{eq:BVP_cost}
\end{align}

\noindent where the subscripts $0$ and $f$ represent the initial and final values, respectively, and the transversality conditions are a result of free final mass and true longitude. The initial guess for the costates of the averaged trajectory was all zeros and the initial guess for the costates of the unaveraged trajectory was the optimal initial costates of the averaged trajectory. The optimized initial costate values for both trajectories are given in the Appendix.

\begin{table}[!t]
    \begin{minipage}{.5\linewidth}
      \footnotesize
      \centering
      \captionof{table}{GTO to GEO parameters.}
        \begin{tabular}{lcc}
        \hhline{===}
         Parameter & Value  \\
         \hline
        Semimajor axis \{GTO, GEO\} & \{24505 km, 42165 km\} \\
        Eccentricity \{GTO, GEO\} & \{0.725, 0\} \\
        Inclination \{GTO, GEO\} & \{28.5 deg, 0 deg\} \\
        AOP \{GTO, GEO\} & \{0 deg, 0 deg\} \\
        RAAN \{GTO, GEO\} & \{0 deg, 0 deg\} \\
        TA \{GTO, GEO\} & \{0 deg, free\} \\
        Initial Mass & 100 kg \\
        Isp & 3100 sec \\
        $J_2$ & 0.00108263 \\
        GM & 398600 km$^3$/sec$^2$ \\
        Integration tolerance & $10^{-14}$ \\
        $\{\epsilon_S, \epsilon_E\}$ & $\{3\times10^{-5}, 1\times10^{-5}\}$ \\
        \hhline{===} \\
    \end{tabular}
    \label{tab:GTOtoGEO_BVP}
    \end{minipage}%
    \begin{minipage}{.5\linewidth}
      \footnotesize
      \centering
      \captionof{table}{GTO to GEO transfer details}
        \begin{tabular}{lccc}
        \hhline{===}
         Parameter & 48 rev  transfer & 486 rev transfer \\
         \hline
        Thrust & [0, 0.2] N & [0, 0.01] N\\
        ToF & 30 Days  & 350 Days\\
        Initial epoch & 260280065 ET & 252417665 ET\\
        Quad Param, $q$ & 6 & 8\\
        Final mass (avg) & 93.645 kg& 91.946 kg\\
        Final mass (unavg) & 93.638 kg & -\\
        $\Delta v$ (avg) & 1.996079 km/s & 2.552701 km/s\\
        $\Delta v$ (unavg) & 1.998351 km/s & -\\
        Eclipse arcs & 34 & 389\\
        Coast arcs & 59 & 377\\
        Integration steps (avg) & 281 & 385 \\
        Integration steps (unavg) & 17,756 & - \\
        \hhline{===} \\
    \end{tabular}
    \label{tab:GTOtoGEO_details}
    \end{minipage} 
\end{table}

\begin{figure}[h!]
\centering
\begin{minipage}{0.55\textwidth}
  \includegraphics[width=1\linewidth]{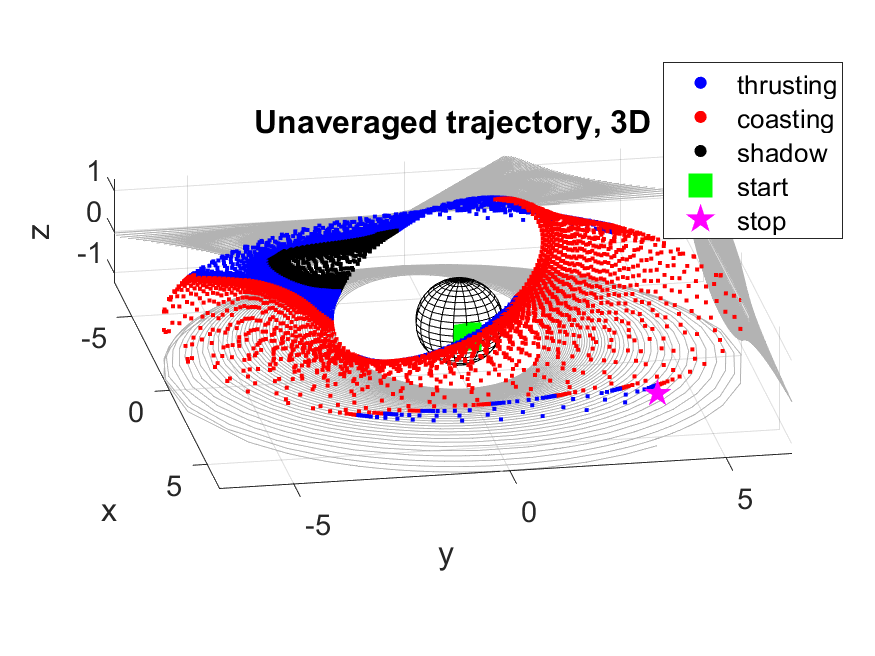}
\end{minipage}
\begin{minipage}{0.49\textwidth}
  \includegraphics[width=1\linewidth]{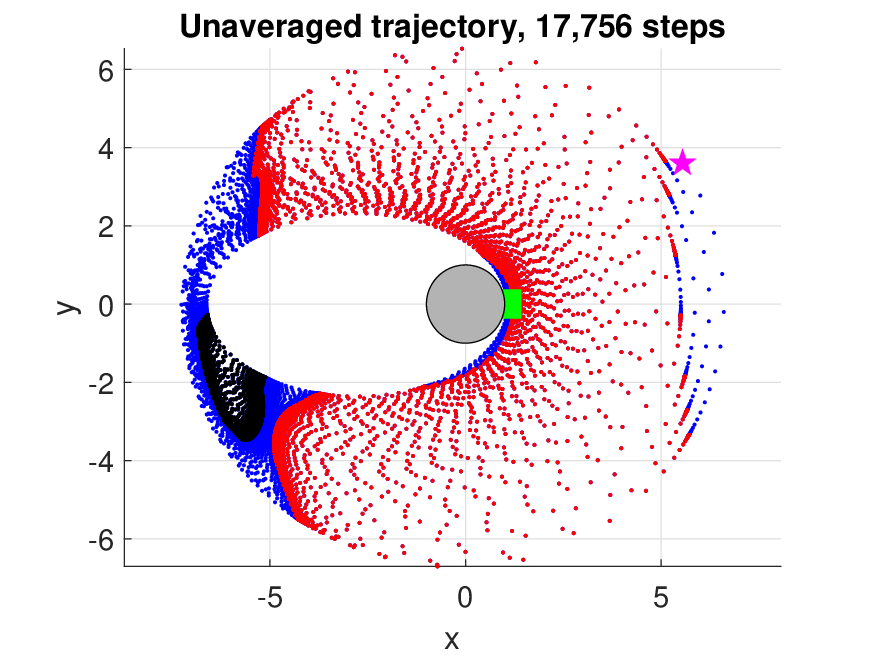}
\end{minipage}
\begin{minipage}{0.49\textwidth}
  \includegraphics[width=1\linewidth]{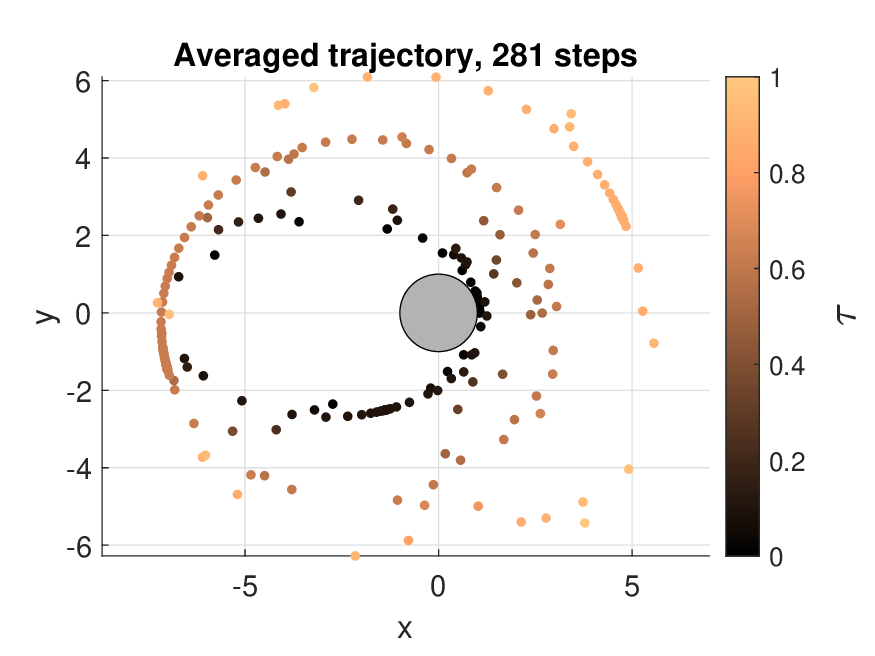}
\end{minipage}
\caption{Trajectory plots for the GTO to GEO transfer. (a: top) 3D, (b: left) unaveraged dynamics, and (c: right) focused averaging dynamics. Units in LU = $6378$ km.}
\label{fig:trajectory}
\end{figure}

\begin{figure}[h!]
\centering
\begin{minipage}{0.48\textwidth}
\centering
  \includegraphics[width=1\linewidth]{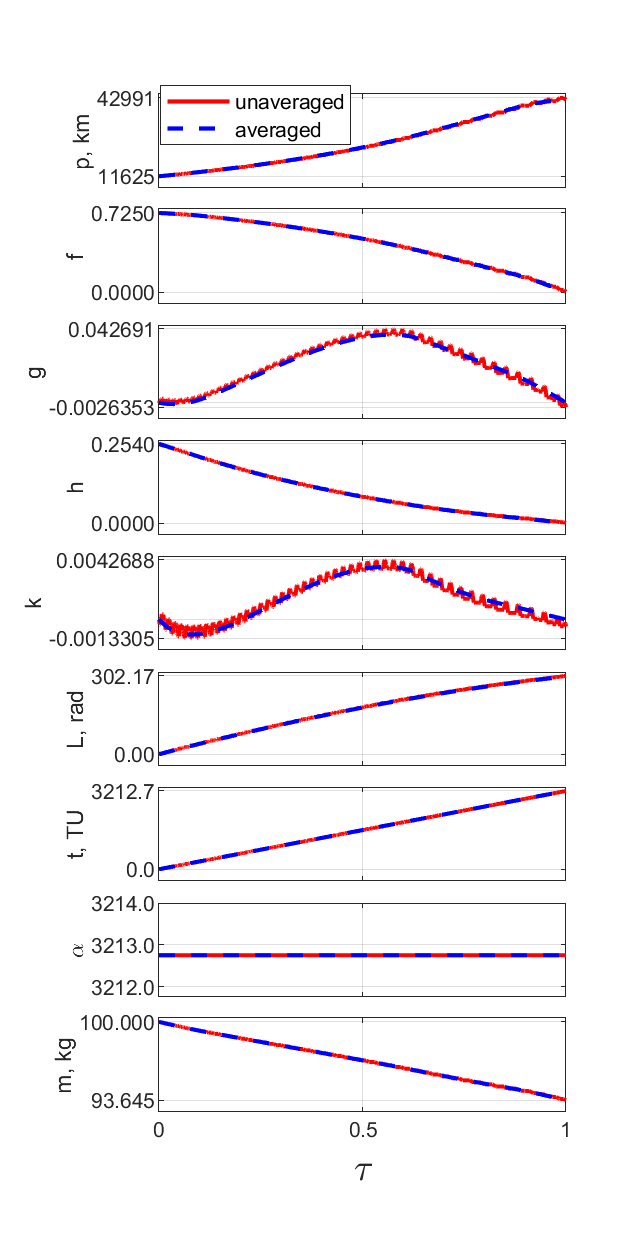}
\end{minipage}
\begin{minipage}{0.48\textwidth}
\centering
  \includegraphics[width=1\linewidth]{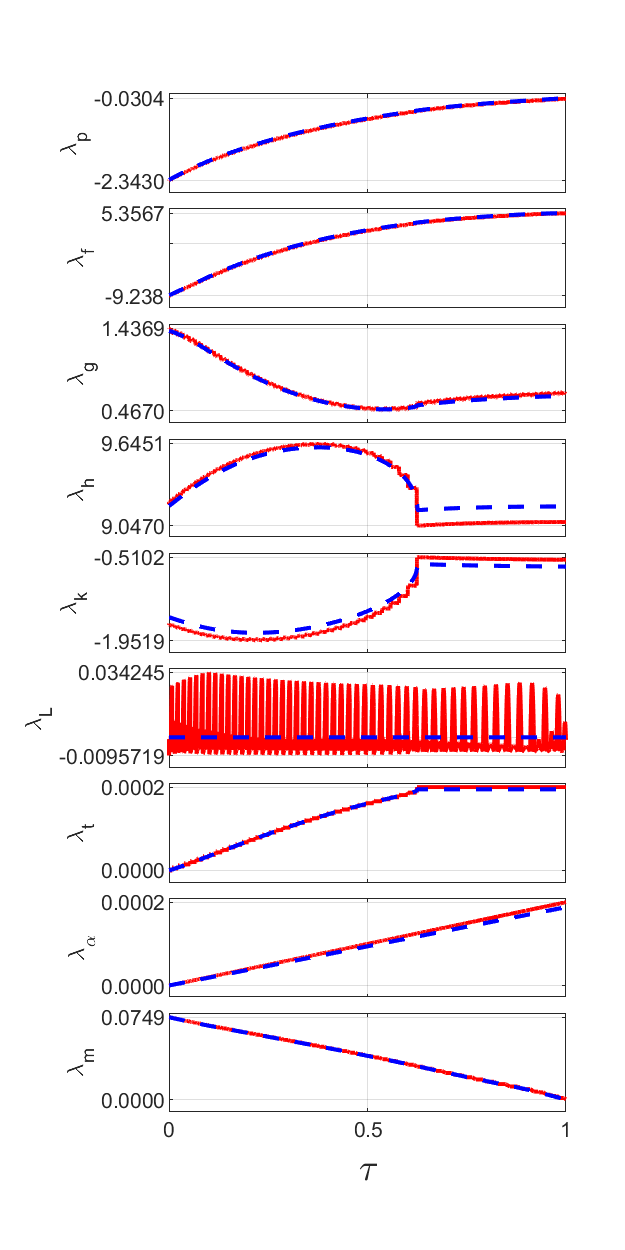}
\end{minipage}
\caption{Comparisons of the (a: left) states and (b: right) costates for the GTO to GEO transfer.}
\label{fig:states_comparison}
\end{figure}

\begin{figure}[h!]
\centering
\begin{minipage}{0.48\textwidth}
\centering
  \includegraphics[width=1\linewidth]{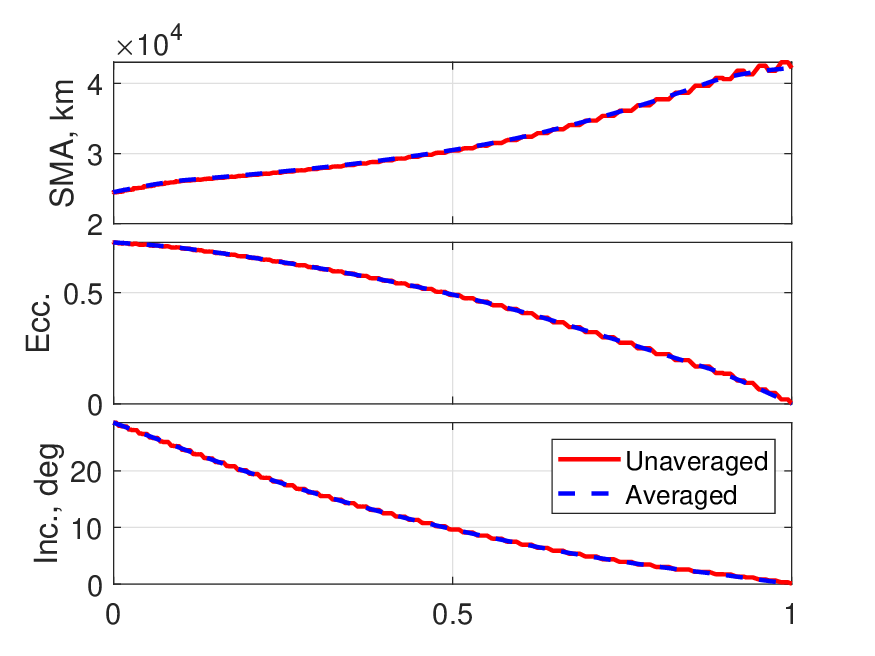}
\end{minipage}
\begin{minipage}{0.48\textwidth}
\centering
  \includegraphics[width=1\linewidth]{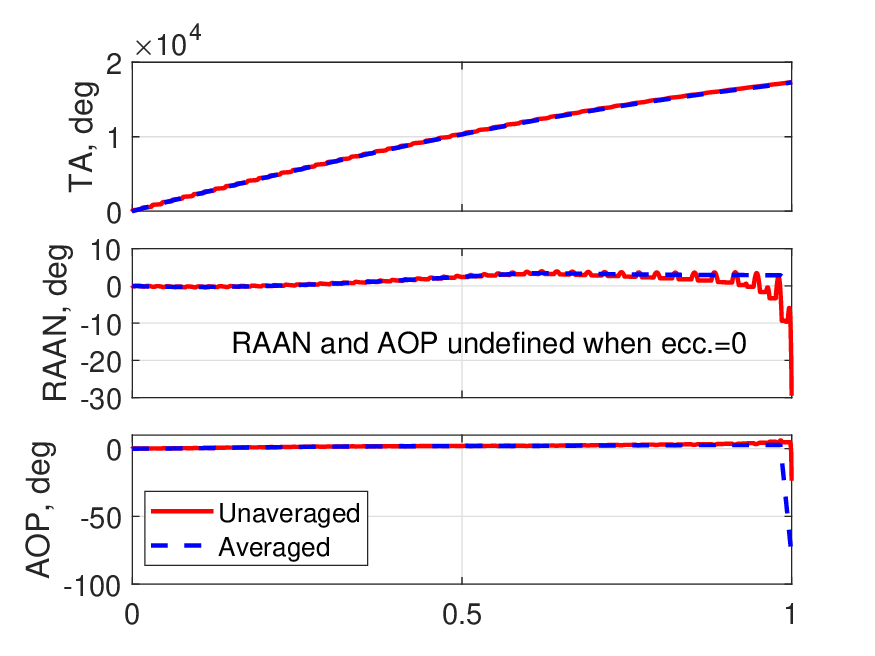}
\end{minipage}
\caption{Comparisons of the classical Keplerian orbital elements for the GTO to GEO transfer.}
\label{fig:keplerian_states}
\end{figure}

\begin{figure}[h!]
\centering
\includegraphics[width=1\linewidth]{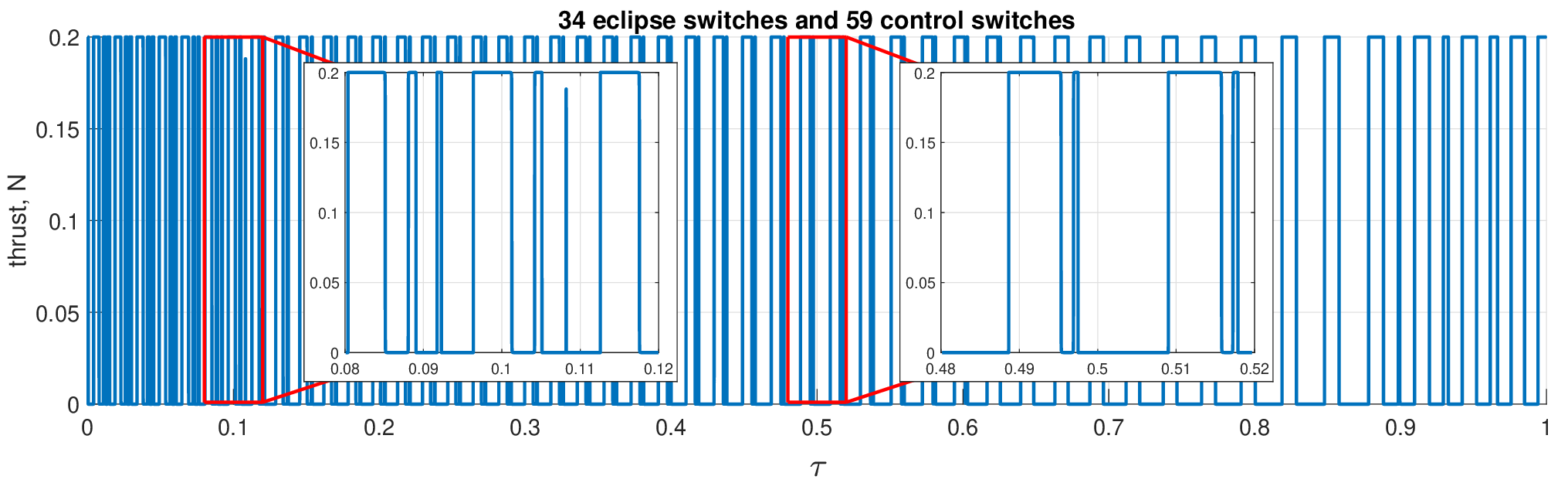}
\caption{Thrust over time for GTO to GEO transfer with unaveraged dynamics.}
\label{fig:thrust}
\end{figure}

The unaveraged trajectory is shown in three dimensions in Fig. \ref{fig:trajectory}(a), with thrusting and shadowing details in Fig. \ref{fig:trajectory}(b). The averaged trajectory is shown with time information in Fig. \ref{fig:trajectory}(c). The sparsity of the points in the averaged trajectory plot show the ability of the variable step integrator to take large steps with averaged dynamics. The variable step integrator takes more frequent steps, though, when the number of thrusting arcs change. The time histories of the states and costates are compared in Fig. \ref{fig:states_comparison}. Additionally, the time histories of states are shown as classical Keplerian orbital elements in Fig. \ref{fig:keplerian_states}. Note that the RAAN and AOP become undefined when eccentricity equals zero. The final mass is a key statistic in the mission design process and a good measure of the ability to use the optimal averaged trajectory as a proxy for the optimal unaveraged trajectory. The optimal unaveraged final mass is 93.638 kilograms (corresponding to 1.998351 km/s $\Delta v$) and the optimal averaged final mass is 93.645 kilograms (corresponding to 1.996079 km/s $\Delta v$). The $\lambda_t$ costate flattens near $\tau=0.65$, because the trajectory no longer encounters any eclipsing and this costate is only affected by the presence of eclipsing. The averaged trajectory closely follows the unaveraged trajectory, providing a validation of the averaged dynamics model as an approximation. Some averaged costates maintain the same time history structure as the unaveraged costates but become slightly out of sync because they have no transversality conditions in solving the boundary value problem.

\begin{figure*}[h!]
    \centering
    \begin{subfigure}[t]{0.49\textwidth}
        \centering
        \includegraphics[width=\textwidth]{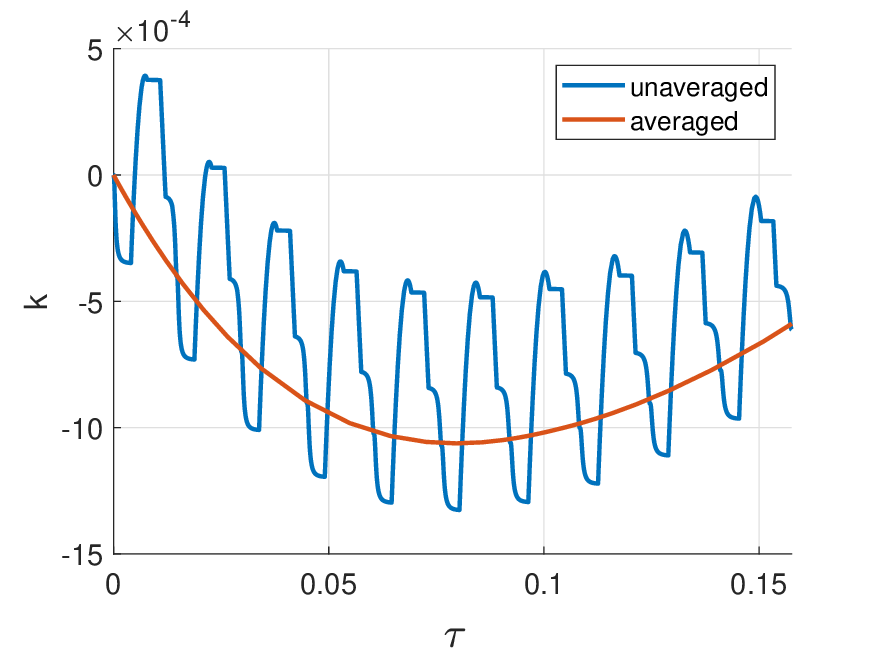}
    \end{subfigure}
    \begin{subfigure}[t]{0.49\textwidth}
        \centering
        \includegraphics[width=\textwidth]{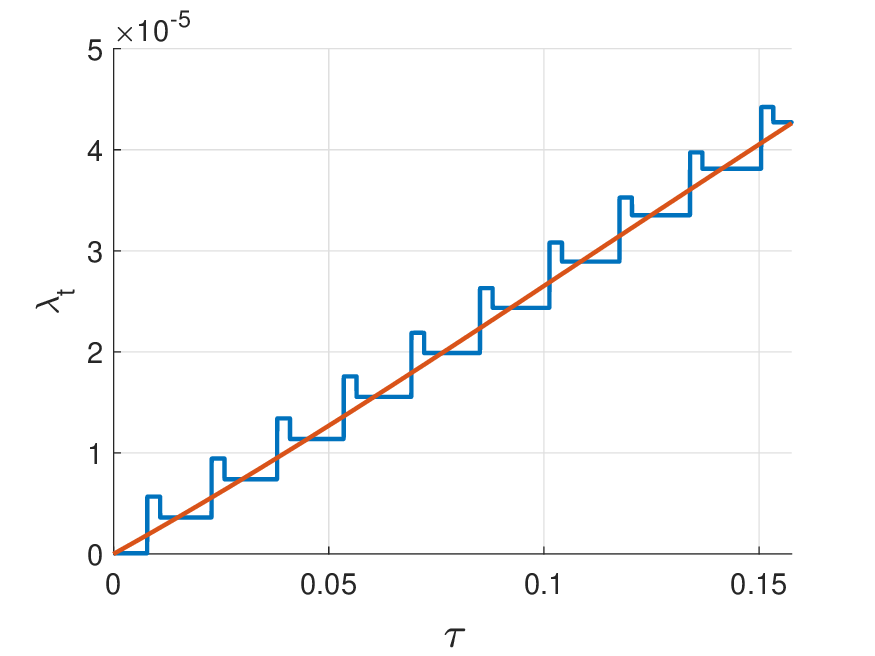}
    \end{subfigure}
    \caption{(a: left) Detailed time history of the first 10 revolutions of state k. (b: right) costate jumps that are smoothly incorporated in the averaged dynamics for $\lambda_t$.}
    \label{avg_jumps}
\end{figure*}

\begin{figure}
    \begin{minipage}{0.48\textwidth}%
        \includegraphics[width=\textwidth]{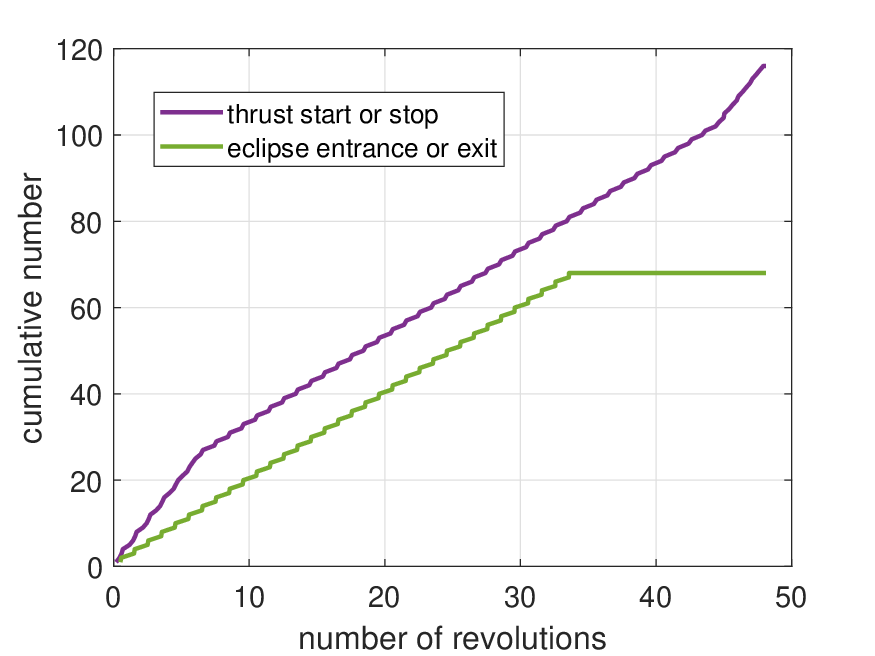}
        \caption{Cumulative number of times trajectory enters or exits an eclipse and starts or stops thrusting due the control.}
        \label{fig:numevents}
    \end{minipage}%
    \hfill
    \begin{minipage}{0.48\textwidth}%
        \includegraphics[width=\textwidth,trim={0cm 7cm 18cm 0cm}]{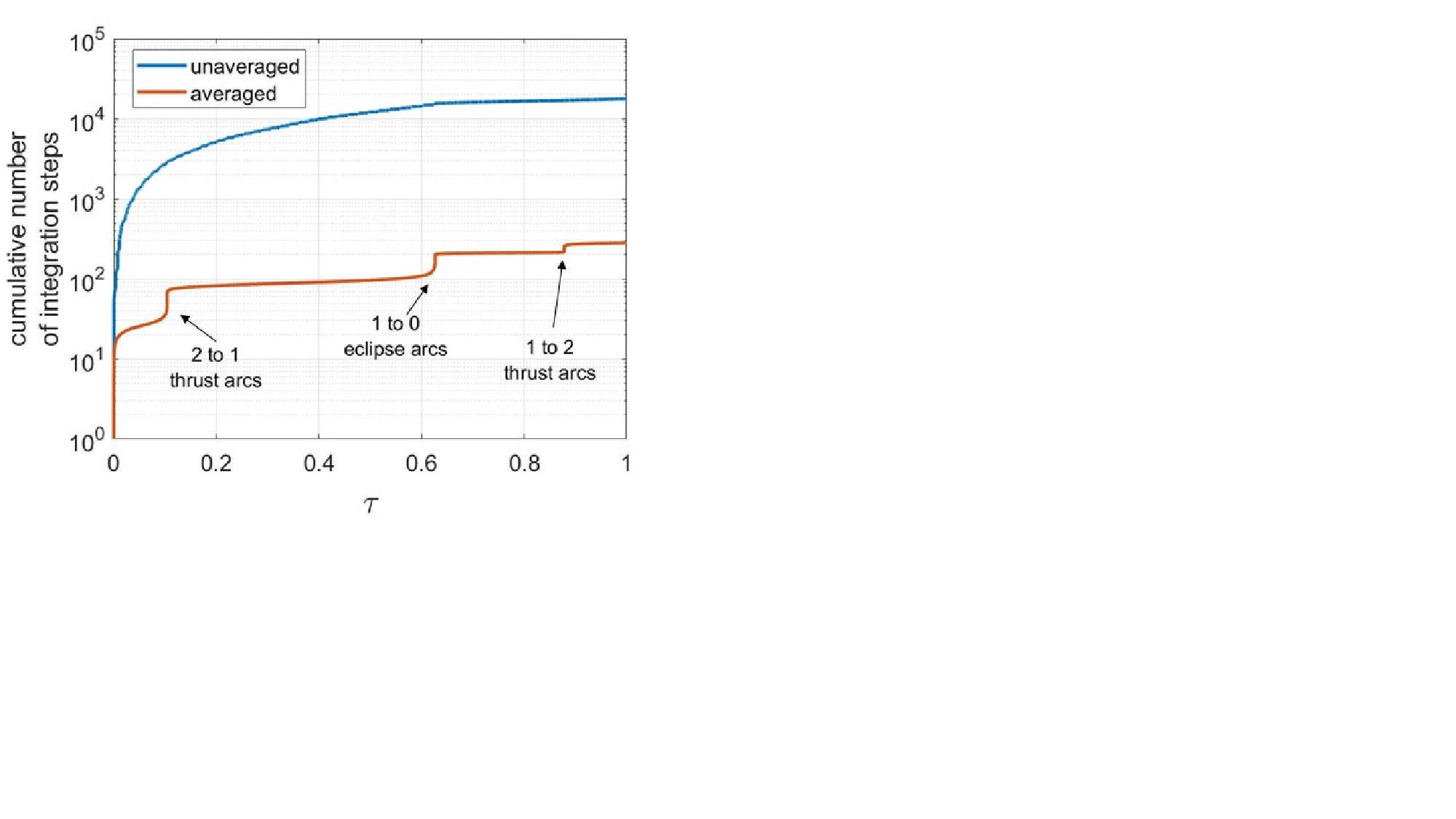}
        \caption{Total of 17,756 steps taken for unaveraged and 281 steps for averaged.}
        \label{fig:numsteps}
    \end{minipage}%
\end{figure}

The thrust profile is shown in Fig. \ref{fig:thrust}. There is one instance during the transfer where the thrust does not completely switch off (highlighted in the first zoom-in panel), because the optimal switch happens so fast that the smoothing parameters are not small enough to capture such a short arc. In Fig. \ref{avg_jumps}(a), more detail is shown for the first 10 revolutions of state $k$ and Fig. \ref{avg_jumps}(b) shows the near-discrete costate jumps of $\lambda_t$ from eclipsing that the averaged dynamics smoothly incorporates. The number of switches and eclipsing events are tracked in Fig. \ref{fig:numevents}. The cumulative number of integration steps are tracked in Fig. \ref{fig:numsteps}, where the unaveraged trajectory requires approximately 60 times more steps. The averaged trajectory takes more frequent steps when the number of thrusting or eclipsing arcs per revolution changes. This slow-down when using a variable step integrator is a inherent downside of averaged dynamics, but typically will only occur a few times during transfers. 

\begin{figure}
    \begin{minipage}{0.48\textwidth}%
        \includegraphics[width=\textwidth]{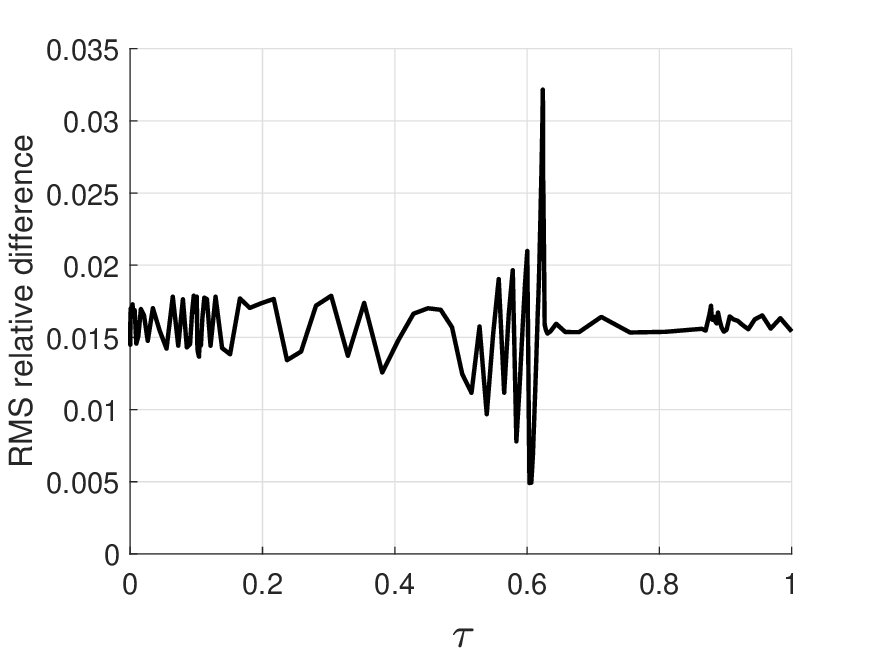}
        \caption{RMS relative difference between averaged and unaveraged trajectory states and costates}
        \label{fig:rms_error}
    \end{minipage}%
    \hfill
    \begin{minipage}{0.48\textwidth}%
        \includegraphics[width=\textwidth]{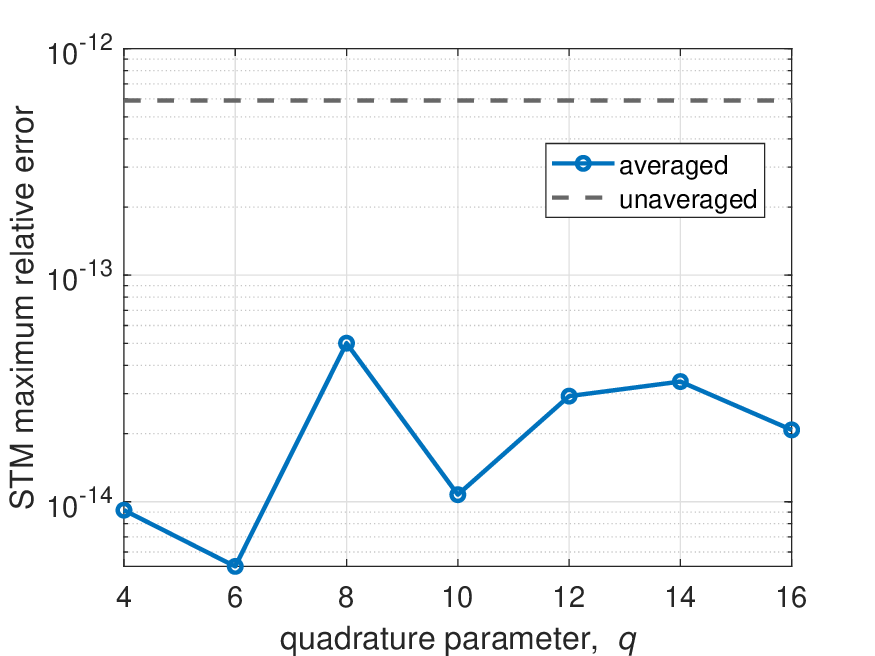}
        \caption{STM accuracy for averaged dynamics compared to complex step derivative.}
        \label{fig:STM_err}
    \end{minipage}%
\end{figure}

The root mean square of the relative difference between the unaveraged and averaged trajectory states and costates is shown in Fig. \ref{fig:rms_error}. This difference is calculated across all states and costates using:
\begin{align}
    \text{difference}(\tau) &= \sqrt{\sum_i^{18} \left[\frac{\text{\textbf{y}}_{\text{avg},i}(\tau) - \text{\textbf{y}}_{\text{unavg},i}(\tau)}{\max(1,\max_{\tau}(|\text{\textbf{y}}_{\text{unavg},i}|))}\right]^2}
\end{align}
This difference is oscillatory and small with a spike when the eclipsing arcs disappear. The difference does not vary significantly with quadrature parameter within reasonable values of $q$ from 6 to 32. 

The variational equations calculated in section \ref{sec:variational} are verified with this trajectory. The STM propagated with the variational equations is compared to the STM calculated using the complex step method, a machine-precision, numerical approach to calculating derivatives described in \cite{martins2003complex} and \cite{pellegrini_STM}. The complex step STM is considered truth for comparison purposes. The relative error for each STM entry is calculated by normalizing by the maximum value in the row or column of each entry:
\begin{align}
    \text{(relative error)}_{i,j} &= \frac{\left|\text{\textbf{A}}_{i,j}-\text{\textbf{B}}_{i,j}\right|}{\max(1,\max_j(\left|\text{\textbf{B}}_i\right|),\max_i(\left|\text{\textbf{B}}_j\right|))}
\end{align}
Extra care is taken to ensure that the trajectories propagated with complex step and with the variational equations follow the exact same path. This path uniformity can be achieved by using fixed step propagation, among other methods \cite{pellegrini_STM}. The maximum relative STM error for different quadrature parameters, $q$, is shown in Fig. \ref{fig:STM_err}. The number of quadrature points per arc does not affect the accuracy of the STM.  Therefore, boundary value problems can be solved efficiently, even with low-fidelity dynamics parameters.


\begin{figure}[h!]
    \begin{minipage}{0.48\textwidth}%
        \includegraphics[width=\textwidth]{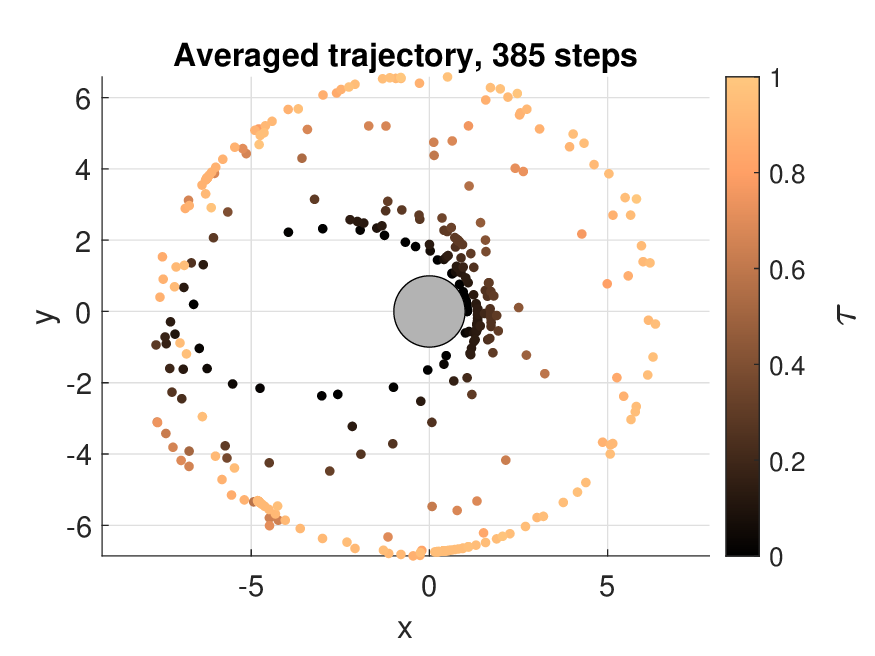}
        \caption{Trajectory plot of the 486 revolution GTO to GEO transfer with averaged dynamics.}
        \label{fig:GTO_long_traj}
    \end{minipage}%
    \hfill
    \begin{minipage}{0.48\textwidth}%
        \includegraphics[width=\textwidth,trim={0cm 7cm 18cm 0cm}]{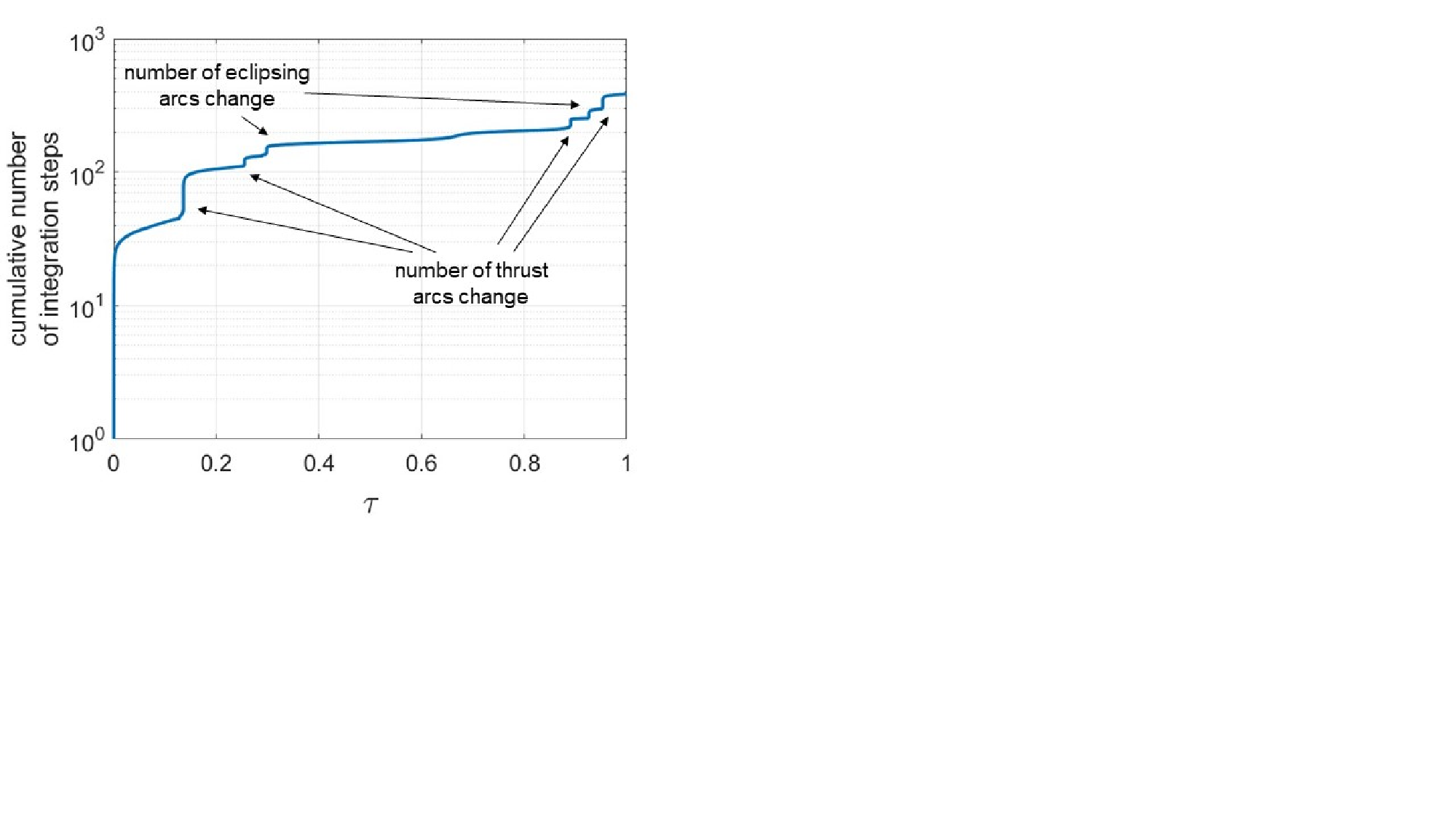}
        \caption{Total number of 385 steps. \bigskip}
        \label{fig:GTO_long_numsteps}
    \end{minipage}%
\end{figure}

\subsection{486-Revolution GTO to GEO example} \label{sec:GTOtoGEO_long}

\begin{figure}[h!]
    \begin{minipage}{0.48\textwidth}%
        \includegraphics[width=\textwidth]{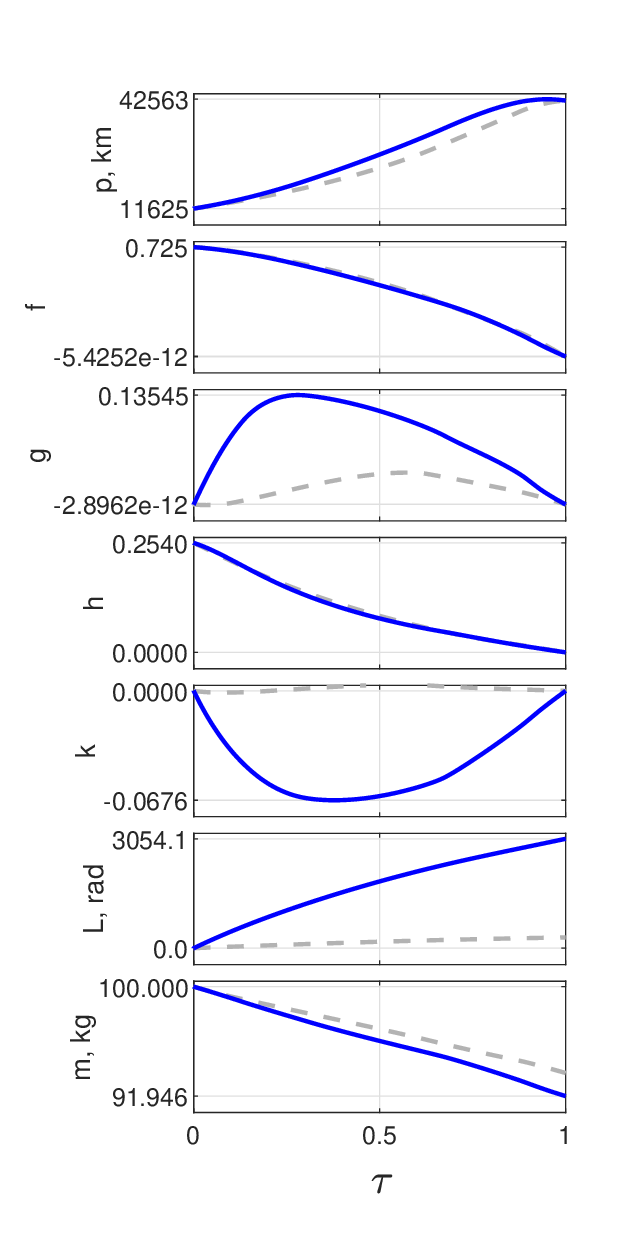}
    \end{minipage}%
    \hfill
    \begin{minipage}{0.48\textwidth}%
        \includegraphics[width=\textwidth]{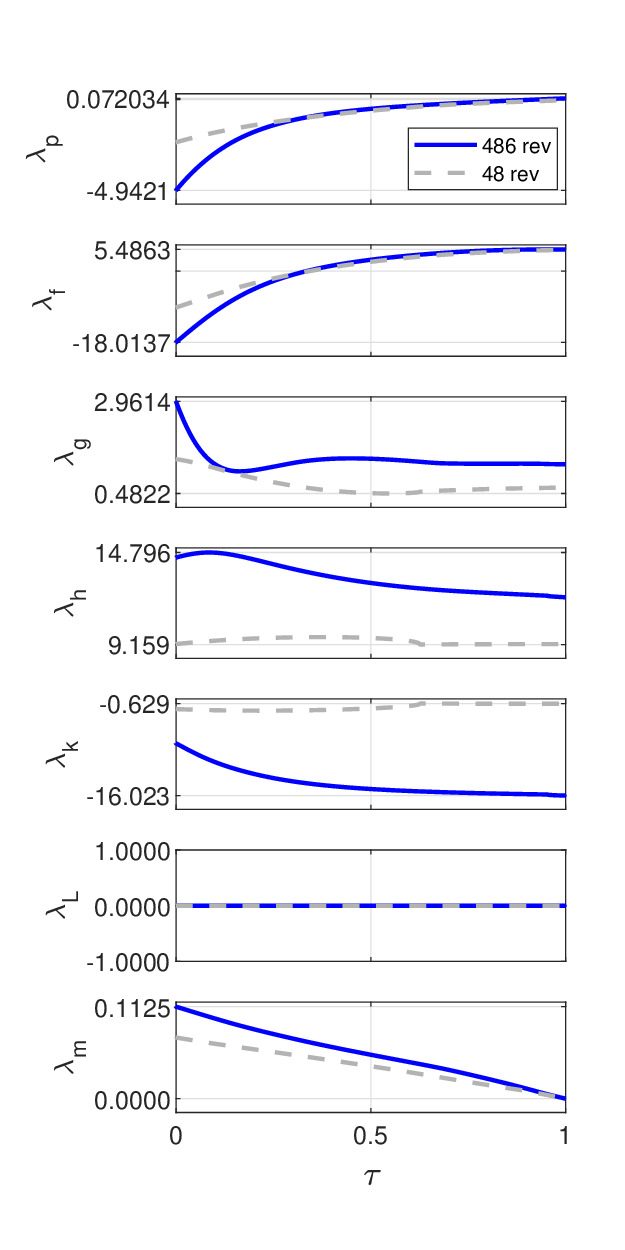}
    \end{minipage}%
    \caption{Plots of the (a: left) states and (b: right) costates for the 486-revolution GTO to GEO transfer with previous example time histories shown in gray.} \label{fig:GTO_long_states}
\end{figure}

\begin{figure}[h!]
\centering
\begin{minipage}{0.48\textwidth}
\centering
  \includegraphics[width=1\linewidth]{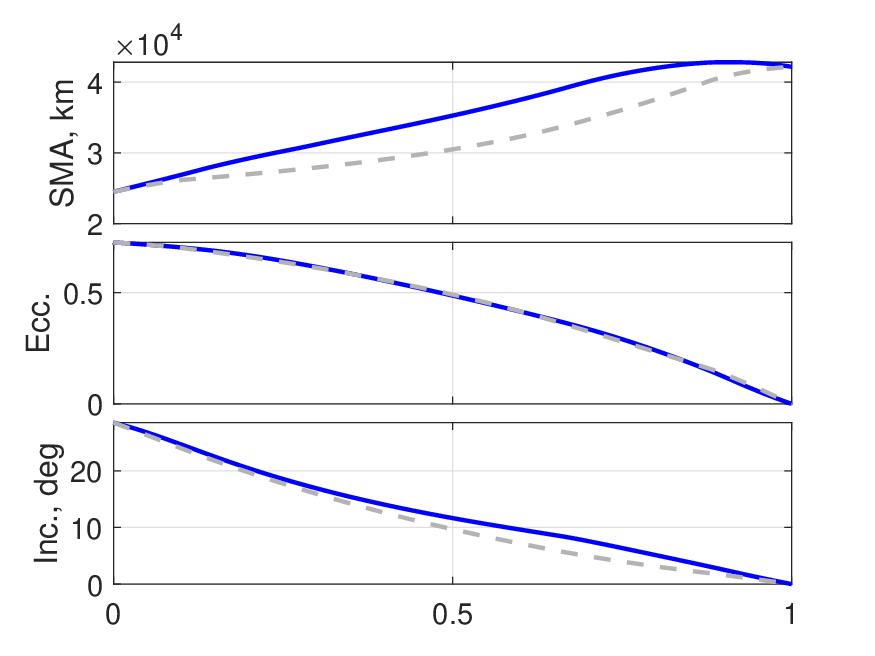}
\end{minipage}
\begin{minipage}{0.48\textwidth}
\centering
  \includegraphics[width=1\linewidth]{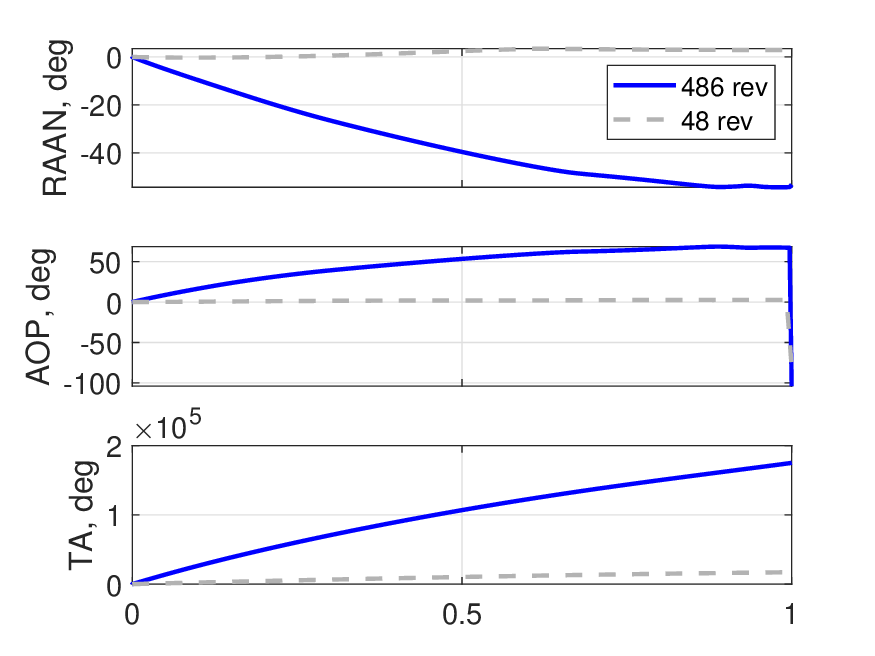}
\end{minipage}
\caption{The classical Keplerian orbital elements for the 486-revolution GTO to GEO transfer.}
\label{fig:keplerian_states_long}
\end{figure}

A 486-revolution, minimum-fuel optimal GTO to GEO transfer is shown to showcase averaged dynamics. The boundary conditions and transfer parameters are identical to Table \ref{tab:GTOtoGEO_BVP} from the previous example. The transfer details are given in Table \ref{tab:GTOtoGEO_details}. The optimized initial costate values are given in the Appendix. The initial guess for the costates of the averaged trajectory was again all zeros. The averaged trajectory is shown with time information in Fig. \ref{fig:GTO_long_traj} and the cumulative number of integration steps are tracked in Fig. \ref{fig:GTO_long_numsteps}. Despite the significant increase in the number of revolutions (486 vs. 48), only about 50 percent more integration steps are required compared to the previous transfer. The averaged dynamics eliminate periodic changes in states and costates, causing the integration steps to scale with qualitative changes like the number of thrusting or eclipsing arcs per revolution instead of the number of total revolutions. The time histories of the states and costates are shown in Fig. \ref{fig:GTO_long_states} with gray dashed lines for the averaged time histories from the previous example. The structure is marginally different from the previous example despite the same boundary conditions, because the eclipsing seasons are different as a result of a different maximum thrust and time-of-flight. Additionally, the time histories of the states are shown as classical Keplerian orbital elements in Fig. \ref{fig:keplerian_states_long} with gray dashed lines for the averaged Keplerian time histories from the previous example. The final mass is 91.946 kg, corresponding to a 2.552701 km/s $\Delta v$.

\section{Conclusion} \label{sec:conclusion}

An improved model is introduced for averaged dynamics and sensitivities for the minimum fuel problem that combines advances from several past works with new developments in the eclipsing constraint, the switching function, an expanded state vector, and variational equations. A singularity in the averaged eclipsing constraint is removed, allowing a variable step integrator to propagate through small eclipsing arcs without stopping, significantly slowing down, or accumulating error. The roots of the switching function are shown to correspond to a sixth-order polynomial in the averaging context, contributing to faster root-finding and thus faster propagation of multi-arc averaged dynamics. This result implies a maximum of three thrusting arcs per revolution in averaged dynamics with deeper implications for the upper bound on number of thrusts per revolution for optimal transfers in general. The variational equations for averaged dynamics are used to calculate an accurate and consistent STM, enabling a smoother optimization problem. These new developments enhance the applicability of the averaged model to optimizing low-thrust, many revolution trajectories with bang-bang control and a discontinuous eclipsing constraint. In practice, the averaged model introduced here can be the final step of a continuation process used in solving the boundary value problem. The initial and intermediary steps of this continuation process can use quadratic cost functions or smooth Heaviside approximations in the averaged model.


\section*{Appendix}

The $\beta$ coefficients for the switching polynomial are shown in Table \ref{tab:beta_coefficients} with parentheses used for algebraic groupings and the initial conditions for the example trajectories are shown in Table \ref{tab:initial_conditions}. The matrices \textbf{A} and \textbf{B} for the equations of motion of modified equinoctial orbital elements are:
\begin{align*}
    \text{\textbf{A}}^\text{T} &= 
    \begin{bmatrix}
        0 & 0 & 0 & 0 & 0 & \sqrt{\mu p} \frac{w^2}{p^2} & 1 & 0
    \end{bmatrix}
\end{align*}
\begin{align*}
    \text{\textbf{B}} &= 
    \begin{bmatrix}
        0 & \frac{2 p}{w} \sqrt{\frac{p}{\mu}}& 0\\
        \sqrt{\frac{p}{\mu}} sin(L) & \sqrt{\frac{p}{\mu}} \frac{1}{w} \{[w+1] \cos(L) + f\} & -\sqrt{\frac{p}{\mu}} \frac{g}{w} [h \sin(L) - k \cos(L)] \\
        -\sqrt{\frac{p}{\mu}} \cos(L) & \sqrt{\frac{p}{\mu}} \frac{1}{w} \{[w+1] \sin(L) + g\} & \sqrt{\frac{p}{\mu}} \frac{f}{w} [h \sin(L) - k \cos(L)] \\
        0 & 0 & \sqrt{\frac{p}{\mu}} \frac{s^2 \cos(L)}{2 w} \\
        0 & 0 & \sqrt{\frac{p}{\mu}} \frac{s^2 \sin(L)}{2 w} \\
        0 & 0 & \sqrt{\frac{p}{\mu}} \frac{1}{w} [h \sin(L) - k \cos(L)] \\
        0 & 0 & 0 \\
        0 & 0 & 0
    \end{bmatrix}
\end{align*}

where $w = 1 + f \cos(L) + g \sin(L)$ and $s^2 = [1 + h^2 + k^2]$.

\begin{table}[h!]
\scriptsize
\centering
\captionof{table}{Beta coefficients for switching polynomial.}
\begin{tabular}{l c}
    \hhline{=|=}
    \begin{normalsize} $\beta_1$ \end{normalsize} & \makecell{$-16 c^{2} p^{3} \lambda_{p}^{2}-32 \lambda_{p} \left(\lambda_{f} f +\frac{1}{2} \lambda_{g} g +\lambda_{f}\right) c^{2} p^{2} -16 \left(\left(\frac{1}{4} \lambda_{g}^{2}+\lambda_{f}^{2}+\frac{1}{4} \lambda_{g}^{2} k^{2}\right) f^{2}+\left(-\frac{\lambda_{g} \lambda_{h} k^{3}}{4}-\frac{\lambda_{g} \left(\lambda_{f} g -\lambda_{L}\right) k^{2}}{2}-\frac{\lambda_{g} \lambda_{h} \left(h^{2}+1\right) k}{4}+\frac{\lambda_{g}^{2}}{2}+2 \lambda_{f}^{2}+\lambda_{f} \lambda_{g} g \right) f \right.$ \\ $\left.+\frac{\lambda_{h}^{2} k^{4}}{16}+\frac{\lambda_{h} \left(\lambda_{f} g -\lambda_{L}\right) k^{3}}{4}+\left(\left(\frac{1}{8}+\frac{h^{2}}{8}\right) \lambda_{h}^{2}+\frac{\left(\lambda_{f} g -\lambda_{L}\right)^{2}}{4}\right) k^{2}+\frac{\lambda_{h} \left(h^{2}+1\right) \left(\lambda_{f} g -\lambda_{L}\right) k}{4}+\frac{\left(h^{2}+1\right)^{2} \lambda_{h}^{2}}{16}+\lambda_{f}^{2}+\lambda_{f} \lambda_{g} g+\frac{\lambda_{g}^{2} \left(g^{2}+1\right)}{4}\right) c^{2} p+4 \mu m^{2} \left(f +1\right)^{2} \left(\lambda_{m}-1\right)^{2}$}\\
    \hhline{--}
    \begin{normalsize} $\beta_2$ \end{normalsize} & \makecell{$32 \lambda_{p} \left(g \lambda_{f}+\lambda_{g} \left(f -2\right)\right) c^{2} p^{2}+32 c^{2} \left(-\frac{\lambda_{f} \left(h k \lambda_{f}-\lambda_{g}\right) g^{2}}{2}+\left(-\frac{\lambda_{f} \lambda_{h} h^{3}}{4}+\frac{\lambda_{f} \lambda_{k} h^{2} k}{4}+\left(-\frac{\lambda_{h} k^{2}}{4}+\left(f \lambda_{g}+\lambda_{L}\right) k -\frac{\lambda_{h}}{4}\right) \lambda_{f} h +\frac{\lambda_{f} \lambda_{k} k^{3}}{4}+\frac{\lambda_{f} \lambda_{k} k}{4}\right.\right.$ \\ $\left.\left. +\left(f -\frac{3}{2}\right) \lambda_{g}^{2}+\lambda_{f}^{2} \left(f -1\right)\right) g +\frac{\lambda_{h} \lambda_{k} h^{4}}{8}+\frac{\lambda_{h} \left(f \lambda_{g}+\lambda_{L}\right) h^{3}}{4}-\frac{\lambda_{k} \left(-\lambda_{h} k^{2}+\left(f \lambda_{g}+\lambda_{L}\right) k -\lambda_{h}\right) h^{2}}{4}-\frac{\left(f \lambda_{g}+\lambda_{L}\right) \left(-\frac{\lambda_{h} k^{2}}{2}+\left(f \lambda_{g}+\lambda_{L}\right) k-\frac{\lambda_{h}}{2}\right) h}{2}\right.$\\$\left.+\frac{\lambda_{h} \lambda_{k} k^{4}}{8}-\frac{\lambda_{k} \left(f \lambda_{g}+\lambda_{L}\right) k^{3}}{4}+\frac{\lambda_{h} \lambda_{k} k^{2}}{4}-\frac{\lambda_{k} \left(f \lambda_{g}+\lambda_{L}\right) k}{4}+\frac{\lambda_{f} \left(f -1\right) \left(f -3\right) \lambda_{g}}{2}+\frac{\lambda_{h} \lambda_{k}}{8}\right) p -16 \mu g m^{2} \left(\lambda_{m}-1\right)^{2} \left(f -1\right)$} \\
    \hhline{--}
    \begin{normalsize} $\beta_3$ \end{normalsize} & \makecell{$-48 c^{2} p^{3} \lambda_{p}^{2}-32 \left(\lambda_{f} f +\frac{7}{2} \lambda_{g} g -\lambda_{f}\right) \lambda_{p} c^{2} p^{2}-16 \left(\left(-\frac{\lambda_{h}^{2}}{16}+\frac{\lambda_{k}^{2}}{4}\right) h^{4}-\lambda_{k} \left(-f \lambda_{g}+g \lambda_{f}-\lambda_{L}\right) h^{3}+\left(\left(\frac{\lambda_{k}^{2}}{2}-\frac{\lambda_{h}^{2}}{8}\right) k^{2}-\frac{\lambda_{h} \left(-f \lambda_{g}+g \lambda_{f}-\lambda_{L}\right) k}{4}+\lambda_{g}^{2} f^{2}\right.\right.$\\$\left.\left.-2 \lambda_{g} \left(g \lambda_{f}-\lambda_{L}\right) f +\lambda_{f}^{2} g^{2}-2 \lambda_{f} \lambda_{L} g -\frac{\lambda_{h}^{2}}{8}+\frac{\lambda_{k}^{2}}{2}+\lambda_{L}^{2}\right) h^{2}-\lambda_{k} \left(k^{2}+1\right) \left(-f \lambda_{g}+g \lambda_{f}-\lambda_{L}\right) h +\left(-\frac{\lambda_{h}^{2}}{16}+\frac{\lambda_{k}^{2}}{4}\right) k^{4}\right.$\\$\left.-\frac{\lambda_{h} \left(-f \lambda_{g}+g \lambda_{f}-\lambda_{L}\right) k^{3}}{4}+\left(-\frac{\lambda_{g}^{2} f^{2}}{4}+\frac{\lambda_{g} \left(g \lambda_{f}-\lambda_{L}\right) f}{2}+\frac{\lambda_{f} \lambda_{L} g}{2}-\frac{\lambda_{h}^{2}}{8}+\frac{\lambda_{k}^{2}}{2}-\frac{\lambda_{L}^{2}}{4}-\frac{\lambda_{f}^{2} g^{2}}{4}\right) k^{2}-\frac{\lambda_{h} \left(-f \lambda_{g}+g \lambda_{f}-\lambda_{L}\right) k}{4}-\frac{\lambda_{g}^{2} f^{2}}{4}+\left(-\frac{5}{2} \lambda_{g}^{2}+3 \lambda_{f} \lambda_{g} g \right) f \right.$\\$\left.+\left(\lambda_{f}^{2}+\frac{15 \lambda_{g}^{2}}{4}\right) g^{2}-5 \lambda_{f} \lambda_{g} g +\frac{\lambda_{k}^{2}}{4}+\frac{15 \lambda_{g}^{2}}{4}-\frac{\lambda_{h}^{2}}{16}\right) c^{2} p-4 \mu m^{2} \left(f^{2}-4 g^{2}+2 f -3\right) \left(\lambda_{m}-1\right)^{2}$}\\
    \hhline{--}
    \begin{normalsize} $\beta_4$ \end{normalsize} & \makecell{$\left(32 \mu \left(\lambda_{m}-1\right)^{2} m^{2}-160 c^{2} p \lambda_{g}^{2}\right) g -128 c^{2} p^{2} \lambda_{p} \lambda_{g}$} \\    
    \hhline{--}
    \begin{normalsize} $\beta_5$ \end{normalsize} & \makecell{$-48 c^{2} p^{3} \lambda_{p}^{2}-32 \left(\lambda_{f} f +\frac{7}{2} \lambda_{g} g +\lambda_{f}\right) \lambda_{p} c^{2} p^{2}-16 \left(\left(-\frac{\lambda_{h}^{2}}{16}+\frac{\lambda_{k}^{2}}{4}\right) h^{4}-\lambda_{k} \left(-f \lambda_{g}+g \lambda_{f}-\lambda_{L}\right) h^{3}+\left(\left(\frac{\lambda_{k}^{2}}{2}-\frac{\lambda_{h}^{2}}{8}\right) k^{2}-\frac{\lambda_{h} \left(-f \lambda_{g}+g \lambda_{f}-\lambda_{L}\right) k}{4}+\lambda_{g}^{2} f^{2}\right.\right.$\\$\left.\left.-2 \lambda_{g} \left(g \lambda_{f}-\lambda_{L}\right) f +\lambda_{f}^{2} g^{2}-2 \lambda_{f} \lambda_{L} g -\frac{\lambda_{h}^{2}}{8}+\frac{\lambda_{k}^{2}}{2}+\lambda_{L}^{2}\right) h^{2}-\lambda_{k} \left(k^{2}+1\right) \left(-f \lambda_{g}+g \lambda_{f}-\lambda_{L}\right) h +\left(-\frac{\lambda_{h}^{2}}{16}+\frac{\lambda_{k}^{2}}{4}\right) k^{4}-\frac{\lambda_{h} \left(-f \lambda_{g}+g \lambda_{f}-\lambda_{L}\right) k^{3}}{4}\right.$\\$\left.+\left(-\frac{\lambda_{g}^{2} f^{2}}{4}+\frac{\lambda_{g} \left(g \lambda_{f}-\lambda_{L}\right) f}{2}+\frac{\lambda_{f} \lambda_{L} g}{2}-\frac{\lambda_{h}^{2}}{8}+\frac{\lambda_{k}^{2}}{2}-\frac{\lambda_{L}^{2}}{4}-\frac{\lambda_{f}^{2} g^{2}}{4}\right) k^{2}-\frac{\lambda_{h} \left(-f \lambda_{g}+g \lambda_{f}-\lambda_{L}\right) k}{4}-\frac{\lambda_{g}^{2} f^{2}}{4}+\left(\frac{5}{2} \lambda_{g}^{2}+3 \lambda_{f} \lambda_{g} g \right) f +\left(\lambda_{f}^{2}+\frac{15 \lambda_{g}^{2}}{4}\right) g^{2}+5 \lambda_{f} \lambda_{g} g\right.$\\$\left. +\frac{\lambda_{k}^{2}}{4}+\frac{15 \lambda_{g}^{2}}{4}-\frac{\lambda_{h}^{2}}{16}\right) c^{2} p-4 \mu m^{2} \left(f^{2}-4 g^{2}-2 f -3\right) \left(\lambda_{m}-1\right)^{2}$} \\
    \hhline{--}
    \begin{normalsize} $\beta_6$ \end{normalsize} & \makecell{$-32 \lambda_{p} \left(g \lambda_{f}+\lambda_{g} \left(f +2\right)\right) c^{2} p^{2}-32 \left(-\frac{\lambda_{f} \left(h k \lambda_{f}-\lambda_{g}\right) g^{2}}{2}+\left(-\frac{\lambda_{f} \lambda_{h} h^{3}}{4}+\frac{\lambda_{f} \lambda_{k} h^{2} k}{4}+\left(-\frac{\lambda_{h} k^{2}}{4}+\left(f \lambda_{g}+\lambda_{L}\right) k -\frac{\lambda_{h}}{4}\right) \lambda_{f} h +\frac{\lambda_{f} \lambda_{k} k^{3}}{4}+\frac{\lambda_{f} \lambda_{k} k}{4}\right.\right.$\\$\left.\left.+\left(f +\frac{3}{2}\right) \lambda_{g}^{2}+\lambda_{f}^{2} \left(f +1\right)\right) g+\frac{\lambda_{h} \lambda_{k} h^{4}}{8}+\frac{\lambda_{h} \left(f \lambda_{g}+\lambda_{L}\right) h^{3}}{4}-\frac{\lambda_{k} \left(-\lambda_{h} k^{2}+\left(f \lambda_{g}+\lambda_{L}\right) k -\lambda_{h}\right) h^{2}}{4}-\frac{\left(f \lambda_{g}+\lambda_{L}\right) \left(-\frac{\lambda_{h} k^{2}}{2}+\left(f \lambda_{g}+\lambda_{L}\right) k -\frac{\lambda_{h}}{2}\right) h}{2}+\frac{\lambda_{h} \lambda_{k} k^{4}}{8}\right.$\\$\left.-\frac{\lambda_{k} \left(f \lambda_{g}+\lambda_{L}\right) k^{3}}{4}+\frac{\lambda_{h} \lambda_{k} k^{2}}{4}-\frac{\lambda_{k} \left(f \lambda_{g}+\lambda_{L}\right) k}{4}+\frac{\lambda_{f} \left(f +3\right) \left(f +1\right) \lambda_{g}}{2}+\frac{\lambda_{h} \lambda_{k}}{8}\right) c^{2} p +16 \mu g m^{2} \left(\lambda_{m}-1\right)^{2} \left(f +1\right)$}\\
    \hhline{--}
    \begin{normalsize} $\beta_7$ \end{normalsize} & \makecell{$-16 c^{2} p^{3} \lambda_{p}^{2}-32 \lambda_{p} \left(\lambda_{f} f +\frac{1}{2} \lambda_{g} g +\lambda_{f}\right) c^{2} p^{2}-16 \left(\left(\frac{1}{4} \lambda_{g}^{2}+\lambda_{f}^{2}+\frac{1}{4} \lambda_{g}^{2} k^{2}\right) f^{2}+\left(-\frac{\lambda_{g} \lambda_{h} k^{3}}{4}-\frac{\lambda_{g} \left(g \lambda_{f}-\lambda_{L}\right) k^{2}}{2}-\frac{\lambda_{g} \lambda_{h} \left(h^{2}+1\right) k}{4}+\frac{\lambda_{g}^{2}}{2}+2 \lambda_{f}^{2}+\lambda_{f} \lambda_{g} g \right) f\right.$\\$\left. +\frac{\lambda_{h}^{2} k^{4}}{16}+\frac{\lambda_{h} \left(g \lambda_{f}-\lambda_{L}\right) k^{3}}{4}+\left(\left(\frac{1}{8}+\frac{h^{2}}{8}\right) \lambda_{h}^{2}+\frac{\left(g \lambda_{f}-\lambda_{L}\right)^{2}}{4}\right) k^{2}+\frac{\lambda_{h} \left(h^{2}+1\right) \left(g \lambda_{f}-\lambda_{L}\right) k}{4}+\frac{\left(h^{2}+1\right)^{2} \lambda_{h}^{2}}{16}+\lambda_{f}^{2}+\lambda_{f} \lambda_{g} g+\frac{\lambda_{g}^{2} \left(g^{2}+1\right)}{4}\right) c^{2} p+4 \mu m^{2} \left(f +1\right)^{2} \left(\lambda_{m}-1\right)^{2}$}\\
    \hhline{=|=} \\
\end{tabular} 
\label{tab:beta_coefficients}
\end{table}

\begin{table}[h!]
      \footnotesize
      \centering
      \captionof{table}{Initial conditions for examples.}
        \begin{tabular}{lcccc}
        \hhline{====}
         State Variable & $\text{\textbf{y}}_{0,\text{unavg}}$ (48 rev) & $\text{\textbf{y}}_{0,\text{avg}}$ (48 rev)& $\text{\textbf{y}}_{0,\text{avg}}$ (486 rev)\\
         \hline
        $p$ & 1.822602598777046 DU & 1.822602598777046 DU & 1.822602598777046 DU\\
        $f$ & 7.25 & 7.25 & 7.25 \\
        $g$ & 0 & 0 & 0 \\
        $h$ & 0.253967646474944 & 0.253967646474944 & 0.253967646474944 \\
        $k$ & 0 & 0 & 0 \\
        $L$ & 0 & 0 & 0 \\
        $t$ & 0 TU & 0 TU & 0 TU \\
        $\alpha$ & 3212.749578552824 TU & 3212.749578552824 TU & 37482.07841644961 TU \\
        $m$ & 100 kg & 100 kg & 100 kg \\
        $\uplambda_p$ & -2.343040603625385 & -2.321725879137949 & -4.942062461563569 \\
        $\uplambda_f$ & -9.220118459700823 & -9.199452707456160 & -18.013685209796851 \\
        $\uplambda_g$ & 1.436912144559582 & 1.406360623157848 & 2.961361665662714 \\
        $\uplambda_h$ & 9.216182583479107 & 9.188890978432537 & 14.486291337614485 \\
        $\uplambda_k$ & -1.673582533704489 & -1.548641252837620 & -7.295171956893319 \\
        $\uplambda_L$ & -0.000976581184152 & 0 & 0 \\
        $\uplambda_t$ & 0.000000059611091 & 0.000000000006312 & -0.000000000842414 \\
        $\uplambda_{\alpha}$ & 0.000000014901161 & 0 & 0 \\
        $\uplambda_m$ & 0.074925287016831 & 0.074834309858591 & 0.112534939505549 \\
        \hhline{====} \\
    \end{tabular}
    \label{tab:initial_conditions}
\end{table}

\section*{Funding Sources}

This work was supported in part by a NASA Space Technology Graduate Research Opportunity, grant number 80NSSC22K1209.

\section*{Acknowledgments}
The authors thank Steve McCarty for technical motivation and advice.

\pagebreak

\bibliography{sample}

\end{document}